\documentclass[final,leqno]{siamltex704}
\usepackage{epsfig}
\usepackage{amsmath}
\usepackage{amssymb}
\usepackage{graphicx}
\usepackage[notcite,notref]{showkeys}

\newtheorem{algorithm}{Weak Galerkin Algorithm}
\newcommand{\bq}{\begin{equation}}
\newcommand{\eq}{\end{equation}}

\newcommand{\bt}{{\bf t}}

\def\T{{\mathcal T}}

\def\bn{{\bf n}}
\def\bq{{\bf q}}

\def\3bar{{|\hspace{-.02in}|\hspace{-.02in}|}}

\title{A Numerical study on the Weak Galerkin method for the Helmholtz equation
with large wave numbers}
\author{lin Mu\thanks{Department of
Applied Science, University of Arkansas at Little Rock,
Little Rock, AR 72204 (lxmu@ualr.edu).} \and Junping
Wang\thanks{Division of Mathematical Sciences, National Science
Foundation, Arlington, VA 22230 (jwang@\break nsf.gov). The research
of Wang was supported by the NSF IR/D program, while working at the
Foundation. However, any opinion, finding, and conclusions or
recommendations expressed in this material are those of the author
and do not necessarily reflect the views of the National Science
Foundation.} \and Xiu Ye\thanks{Department of Mathematics and
Statistics, University of Arkansas at Little Rock, Little Rock, AR
72204 (xxye@ualr.edu). This research of Ye was supported in part by
National Science Foundation Grant DMS-1115097.} \and Shan
Zhao\thanks{Department of Mathematics, University of Alabama,
Tuscaloosa, AL 35487 (szhao@bama.ua.edu). The research of Zhao was
supported in part by National Science Foundation Grant DMS-1016579.
}}
\begin{document}
\maketitle

\begin{abstract}
Weak Galerkin (WG) refers to general finite element methods for
partial differential equations in which differential operators are
approximated by weak forms through the usual integration by parts.
In particular, WG methods allow the use of discontinuous finite
element functions in the algorithm design. One of such examples was
recently introduced in \cite{wy} for solving second order elliptic
problems. The goal of this paper is to apply the WG method of
\cite{wy} to the Helmholtz equation with high wave numbers. Several
test scenarios are designed for a numerical investigation on the
accuracy, convergence, and robustness of the WG method in both
inhomogeneous and homogeneous media over convex and non-convex
domains. Our numerical experiments indicate that weak Galerkin is a
finite element technique that is easy to implement, and provides
very accurate and robust numerical solutions for the Helmholtz
problem with high wave numbers.
\end{abstract}


\pagestyle{myheadings}

\section{Introduction}
In this paper, we explore the use of a weak Galerkin (WG) finite
element method for solving the nonhomogeneous Helmholtz equation
with high wave numbers
\begin{equation}\label{Helmholtz}
-\nabla \cdot (d \nabla u) -k^2u = f, \quad
%\mbox{in}\;\Omega,
\end{equation}
where $k$ is the wave number, $f$ represents a harmonic source, and
$d=d(x,y)$ is a spatial function describing the dielectric
properties of the medium. The Helmholtz equation (\ref{Helmholtz})
governs many macroscopic wave phenomena in the frequency domain
including wave propagation, guiding, radiation and scattering, where
the time-harmonic behavior can be assumed. The numerical solution to the
Helmholtz equation plays a vital role in a wide range of
applications in electromagnetics, optics, and acoustics, such as
antenna analysis and synthesis, radar cross section calculation,
simulation of ground or surface penetrating radar, design of
optoelectronic devices, acoustic noise control, and seismic wave
propagation. However, it remains a challenge to design robust and
efficient numerical algorithms for the Helmholtz equation,
especially when high wave numbers or highly oscillatory solutions
are involved \cite{Zienkiewicz}.

Physically, the Helmholtz problem is usually defined on an unbounded
exterior domain with the so-called Sommerfeld radiation condition
holding at infinity \cite{harari,shen07}
\begin{equation}
\frac{\partial u}{\partial r} - i k u = o \left( r^{\frac{1-m}{2}}
\right), \quad \mbox{as}~ r \to \infty, \; m=2,3.
\end{equation}
where $i=\sqrt{-1}$ is the imaginary unit and $r$ is the radial
direction. Here we have assumed that $d=1$ in far field. In
computational electromagnetics, the exterior domain problems are
often solved numerically by introducing a bounded domain $\Omega$
with an artificial boundary $\partial\Omega$ and imposing certain
boundary conditions on $\partial\Omega$ so that nonphysical
reflections from the boundary can be eliminated or minimized. For
the Helmholtz exterior problems, the non-reflecting condition is
commonly chosen as a Dirichlet-to-Neumann (DtN) mapping, which
relates the wave solution to its derivatives on $\partial\Omega$
\cite{harari,shen07}
\begin{equation}\label{DtN}
d\nabla u \cdot \bn - T(u) = g, \quad \mbox{on}\;
\partial\Omega,
\end{equation}
where $\bn$ denotes the outward normal direction of $\partial\Omega$,
$T$ is a DtN integral operator, and $g$ is a given data function.
For sufficiently large $\Omega$, the nonlocal boundary condition (\ref{DtN})
can be approximated by a Robin boundary condition  \cite{harari,shen07}
\begin{equation}
d\nabla u \cdot \bn - i k u = g,
\quad \mbox{on}\; \partial\Omega,
\end{equation}
which is essentially a first order absorbing boundary condition.

In the present paper, we consider the following prototype Helmholtz
problem
\begin{eqnarray}
-\nabla \cdot (d \nabla u)-k^2u &=& f\quad
\mbox{in}\;\Omega,\label{pde}\\
d\nabla u \cdot \bn- iku&=&g\quad \mbox{on}\;
\partial\Omega.\label{bc}
\end{eqnarray}
Finite element methods for such Helmholtz problems can be classified
as two categories. The first category consists of methods that use
continuous functions to approximate the solution $u$ and the other
refers to methods with discontinuous approximation functions.

Continuous Galerkin (CG) finite element methods employ continuous
piecewise polynomials to approximate the true solution of
(\ref{pde})-(\ref{bc}) and lead to a simple formulation: find $u_h\in
V_h\subset H^1(\Omega)$ satisfying
\begin{equation}\label{cg}
(d\nabla u_h,\nabla
v_h)-k^2(u_h,\;v_h)+ik(u_h,\;v_h)_{\partial\Omega}= (f, v_h)
+(g,v_h)_{\partial\Omega}
\end{equation}
for all $v_h\in V_h$, where $V_h$ is an properly defined finite
element space consisting of continuous piecewise polynomials.

Discontinuous Galerkin (DG) finite element methods for the Helmholtz
equation with $d=1$ (e.g. \cite{fw}) seek $u_h\in V_h\subset
L^2(\Omega)$ satisfying
\begin{eqnarray}
\sum_T(\nabla u_h,\nabla v_h)_T&&-\sum_e ((\{\nabla u_h\}, [v_h])_e+\sigma (\{\nabla v_h\}, [u_h])_e)
-k^2(u_h,\;v_h)\nonumber\\
&&+ik(u_h,\;v_h)_{\partial\Omega}+i(\sum_e\alpha_1h_e^{-1}([u_h],\;[v_h])_e+\sum_e\alpha_2h_e([\frac{\partial u_h}{\partial\bn}],\;[\frac{\partial v_h}{\partial\bn}])_e\nonumber\\
&&+\sum_e\alpha_3h^{-1}_e([\frac{\partial u_h}{\partial\bt}],\;[\frac{\partial v_h}{\partial\bt}])_e)=(f,v_h)+(g,v_h)_{\partial\Omega},\label{dg}
\end{eqnarray}
for all $v_h\in V_h$, where $\alpha_i$, $i$=1,2,3 are penalty
parameters and $V_h$ is a space of discontinuous piecewise
polynomials.

The continuous Galerkin finite element formulation (\ref{cg}) is
natural and easy to implement. However, the errors of continuous
Galerkin finite element solutions deteriorate rapidly when $k$
becomes large. On the other hand, the formulation of DG methods
(\ref{dg}) is complex, and the stability and accuracy of DG methods
heavily depend on the selection of penalty parameters (a good
determination for the parameter values has been an issue for DG
methods).

The objective of the present paper is to introduce a
weak Galerkin (WG) finite element method for solving the
Helmholtz problem with high wave numbers and to investigate the robustness and effectiveness of such WG method
through many carefully designed numerical experiments.
The weak Galerkin finite element formulation was first developed in \cite{wy} for solving the second order elliptic equations. Through rigorous error analysis,
optimal order of convergence of the WG solution in both
discrete $H^1$ norm and $L^2$ norm is established under
minimum regularity assumptions in \cite{wy}.
Nevertheless, no experimental results have been reported in \cite{wy}.

Generally speaking, like the DG method, the WG
finite element method allows one to use discontinuous functions in
the finite element procedure. The main idea behind weak Galerkin
methods lies in the approximation of differential operators by weak
forms for discontinuous finite element functions defined on a
partition of the domain. For example, for the gradient operator, one
may introduce a weak gradient operator $\nabla_d$ over discontinuous
functions $v_h=\{v_0, v_b\}$, where $v_0$ is defined in the interior
of the element and $v_b$ is defined on the boundary of the element.
The weak gradient operator will then be employed to form a weak
Galerkin finite element formulation for (\ref{pde})-(\ref{bc}): find
$u_h\in V_h$ such that for all $v_h\in V_h$ we have
\begin{equation}\label{wg1}
(d\nabla_d u_h,\nabla_d
v_h)-k^2(u_0,\;v_0)+ik(u_b,\;v_b)_{\partial\Omega}= (f, v_0)
+(g,v_b)_{\partial\Omega}.
\end{equation}
The discrete gradient $\nabla_dv$ in (\ref{wg1}) is a variable
locally calculated on each element. The weak Galerkin methods have
many nice features. First, the formulation of WG method is simple,
easy to implement, and involves no penalty parameters for users to
select. The weak Galerkin is equipped with elements of polynomials of
any degree $k\ge 0$. Secondly, the weak Galerkin method conserves
mass with a well defined numerical flux. In some sense, weak
Galerkin finite element methods enjoy both the simplicity of CG
method and the flexibility of DG method.

In the present numerical investigation, we are particularly
interested in the performance of the WG method for solving the Helmholtz
equation with high wave numbers. In general, the numerical
performance of any finite element solution to the Helmholtz equation
depends significantly on the wave number $k$. When $k$ is very large
-- representing a highly oscillatory wave, the mesh size $h$ has to be
sufficiently small for the scheme to resolve the oscillations. To
keep a fixed grid resolution, a natural rule is to choose $kh$ to be
a constant in the mesh refinement, as the wave number $k$ increases
\cite{IhlBabH,bao04}. However, it is known
\cite{IhlBabH,IhlBabHP,BIPS,BabSau} that, even under such a mesh
refinement, the errors of continuous Galerkin finite element
solutions deteriorate rapidly when $k$ becomes larger. This
non-robust behavior with respect to $k$ is known as the ``pollution
effect'' \cite{IhlBabH,IhlBabHP,BIPS,BabSau}. Usually, under the
small magnitude assumption of $kh$, the relative error of a $p$th
order finite element solution in the $H^1$-semi-norm consists of two
parts \cite{IhlBabH,IhlBabHP}, i.e., an error term of the best
approximation behavior like $O(k^p h^p)$ and a pollution error term
behavior like $O(k^{2p+1} h^{2p})$. The error of best approximation
is essentially due to the interpolation error on a discretized grid
and is of bounded magnitude if $kh=\mbox{constant}$. Nevertheless,
the pollution error term dominates when $k$ is large and is
responsible for the non-robustness behavior in the finite element
solutions to the Helmholtz equation.

To the end of eliminating or substantially reducing the pollution
errors, various numerical approaches have been developed in the
literature for solving the Helmholtz equation. For the continuous and
least-squares finite element methods, a popular strategy
\cite{BIPS,BabSau,Melenk,Babuska04,Monk99} to reduce the pollution
error is to include some analytical knowledge of the problem, such
as characteristics of asymptotic or exact solution, in the finite
element space. The resulting local basis functions of non-polynomial
shape yield improved performance
\cite{BIPS,BabSau,Melenk,Babuska04,Monk99}. Similarly, analytical
information is incorporated in the basis functions  of the boundary
element methods to address the high frequency problems
\cite{Giladi,Langdon06,Langdon07}. Based on the geometrical optics
and geometrical theory of diffraction, asymptotic solutions of the
Helmholtz equation that represent important wave propagation
directions are coupled with high frequency oscillations to built
boundary element approximation spaces.
Consequently, the number of degrees of freedom of
such boundary element methods virtually does not depend on the wave
number $k$ for the Helmholtz equation
\cite{Giladi,Langdon06,Langdon07}. In \cite{UWVF03}, plane wave
solutions traveling in a large number of directions have also been
employed as basis functions in the ultra weak variational
formulation (UWVF) \cite{UWVF98}. This allows the use of a coarse
mesh when resolving high oscillatory wave solutions. The UWVF method
can be regarded as an upwinding DG discretization for a first order
system obtained through the introduction of an adjoint field. Based
on such a viewpoint, error estimates of the UWVF on the entire
domain have been recently presented in \cite{Monk08}. Multiple plane
wave basis functions are also employed in a DG method using both low
order elements \cite{Farhat03} and high order elements
\cite{Farhat04}. The weak continuity of the solution at the element
boundaries is enforced by introducing a Lagrange multiplier. Being
more robust than classical variational approaches, this DG method
can resolve high frequency short wave problems very well with a
ratio of up to three wavelengths per element
\cite{Farhat03,Farhat04}.

Besides changing basis functions, there are other improvements that
can be made in the DG framework to reduce the pollution error. By
using piecewise linear polynomials as the basis functions, stable
interior penalty DG methods have been developed in \cite{fw} through
penalizing not only the jumps of the function values, but also the
jumps of the normal and tangential derivatives across the element
edges. Robust results against pollution effect can be achieved
through a careful selection of penalty parameters \cite{fw}. Also
targeting on enhancing stability, DG methods can be established by
formulating the wave equation as a first-order system
\cite{abcm,Chung}. The classical techniques,
such as the enforcement of
weak continuity via fluxes form, can be used in such DG solutions to
the Helmholtz equation. The dispersive and dissipative behavior of
selected DG schemes for both second order Helmholtz equation and the
corresponding first order system have been examined in
\cite{Ainsworth06}. Hybridizable DG method is another efficient
finite element method for the Helmholtz equation \cite{cdg,cgl}. By
introducing new degree of freedoms on the boundary of elements,
parametrization is conducted element by element so that the final
linear system consists of unknowns only from the skeleton of the
mesh. This greatly reduces the size of linear systems compared to
the standard DG scheme. The error analysis \cite{Monk11} shows that
the condition number of the condensed matrix of the hybridizable DG
methods \cite{cdg,cgl} could be independent of the wave number.

The pollution effect can also be alleviated by controlling the
numerical dispersion, i.e., the phase difference between the
numerical and exact waves. This is because the pollution error is
known to be directly related to the dispersion error
\cite{IhlBab,Ainsworth04}, while the best approximation error is
non-dispersive. Thus, the reduction of the dispersive error is
equivalent to the reduction of the pollution error. Consequently,
high order methods, such as local spectral method \cite{bao04},
spectral Galerkin methods \cite{shen05,shen07}, and spectral element
methods \cite{Heikkola,Ainsworth09} are less vulnerable to the
pollution effect, because they produce negligible dispersive errors.

The rest of this paper is organized as follows. In Section 2, we
will introduce a weak Galerkin finite element method for the
Helmholtz equation by following the idea presented in \cite{wy}. In
Section 3, we will present some error estimates for the WG
method for the Helmholtz equation. Finally in Section 4, we shall
present some numerical results obtained from the weak Galerkin
method with various orders.

\section{A Weak Galerkin Finite Element Method}

Let ${\cal T}_h$ be a partition of the domain $\Omega$ with mesh size
$h$. Assume that the partition ${\cal T}_h$ is shape regular so that
the routine inverse inequality in the finite element analysis holds
true (see \cite{ci}). Define $(v,w)_K=\int_K vwdx$ and
$(v,w)_{\partial K}=\int_{\partial K}vwds$.

For each triangle $T\in {\cal T}_h$, let $T^0$ and $\partial T$
denote the interior and boundary of $T$ respectively. Denote by
$P_j(T^0)$ the set of polynomials in $T^0$ with degree no more than
$j$, and $P_\ell(\partial T)$ the set of polynomials on each segment
(edge or face) of $\partial T$ with degree no more than $\ell$. We
emphasize that functions of $P_\ell(\partial T)$ are defined on each
edge/face and there is no continuity required across different
edges/faces. Define a global function $v=\{v_0, v_b\}$ where $v_0$
and $v_b$ represent the values of $v$ on $T^0$ and $\partial T$
respectively for each $T\in {\cal T}_h$. Now we are ready to define
a weak Galerkin finite element space as follows
\begin{equation}\label{vh}
V_h=\left\{ v=\{v_0, v_b\}:\ \{v_0, v_b\}|_{T}\in P_j(T^0)\times P_\ell(\partial T), \forall T\in {\cal T}_h \right\}.
\end{equation}
For each $v=\{v_0,
v_b\}\in V_h$, we define the discrete gradient of $v$ on each
element $T$  by the following equation:
\begin{equation}\label{discrete-weak-gradient-new}
\int_T \nabla_{d} v\cdot q dT = -\int_T v_0 \nabla\cdot q dT+
\int_{\partial T} v_b q\cdot\bn ds,\qquad \forall q\in V_r(T),
\end{equation}
where $V_r(T)$ is a subspace of the set
of vector-valued polynomials of degree no more than $r$ on $T$.

For the purpose of easy demonstration, we consider a special 2-D
weak Galerkin element of $V_h$ and $V_r(T)$ with $j=\ell=k\ge 0$ in
(\ref{vh}) and $V_r(T)=RT_k(T)$. More weak Galerkin elements can be
found in \cite{wy}. Here $RT_k(T)$ is the usual Raviart-Thomas
element \cite{rt} of order $k$ which has the form
\[
RT_k(T)=P_k(T)^2+\tilde{P}_k(T){\bf x},
\]
where $\tilde{P}_k(T)$ is the set of homogeneous polynomials of
degree $k$ and ${\bf x}=(x_1,x_2)$.

\medskip

\begin{algorithm}
A numerical approximation for (\ref{pde})-(\ref{bc}) can be
obtained by seeking $u_h=\{u_0,u_b\}\in V_h$  such that for all $v_h=\{v_0,v_b\}\in V_h$
\begin{equation}\label{WG}
(d\nabla_d u_h,\nabla_d
v_h)-k^2(u_0,\;v_0)+ik(u_b,\;v_b)_{\partial\Omega}= (f, v_0)
+(g,v_b)_{\partial\Omega}.
\end{equation}
\end{algorithm}

In the following, we will use the lowest order weak Galerkin element
($k$=0) as an example to demonstrate how one might implement weak
Galerkin finite element method for solving the Helmholtz problem
(\ref{pde})-(\ref{bc}). Let $N(T)$ and $N(e)$ denote the number of
triangles and the number of edges associated with triangulation
${\cal T}_h$.  Let ${\cal E}_h$ denote the union of the boundaries
of the triangles $T$ of ${\cal T}_h$. Let $\phi_i$ be a function
which takes value one in the interior of triangle $T_i$ and  zero
everywhere else. Let $\psi_j$ be a function that takes value one on
edge $e_j$ and zero everywhere else. The weak Galerkin finite
element space $V_h$ with $k=0$ has the form
\[
V_h={\rm span} \{\phi_1,\cdots,\phi_{N(T)},\psi_1,\cdots,\psi_{N(e)}\}.
\]
The methodology of implementing WG is the same as that for continuous
Galerkin finite element method except that the standard gradient
operator $\nabla$ must be replaced by the discrete gradient operator
$\nabla_d$. The discrete gradient operator is an extension of the
standard gradient operator for smooth functions to non-smooth
functions. Therefore, the key step of implementing weak Galerkin
method is to compute a {\bf local} variable $\nabla_d v$ for $v
=\{v_0,v_b\}\in V_h$ element by element. For the case $k=0$, we have
$\nabla_dv\in RT_0(T)$ on each element $T\in {\cal T}_h$ where
$$
RT_0(T)=\left(\begin{array}{c} a+cx \\b+cy\\\end{array}\right)={\rm span}
\{\theta_1,\theta_1,\theta_3\}.
$$
For example, we can choose $\theta_i$ as follows
\[
\theta_1=\left(\begin{array}{c}1 \\0 \\\end{array}\right), \theta_2=\left(\begin{array}{c}0 \\1 \\\end{array}\right),
\theta_3=\left(\begin{array}{c}x \\y \\\end{array}\right).
\]
Thus on each element $T\in {\cal T}_h$, $\nabla_d v=\sum_{j=1}^3c_j\theta_j$. Using the definition of discrete gradient (\ref{discrete-weak-gradient-new}), we find  $c_j$ by solving the following linear system:
\[
\left(
  \begin{array}{ccc}
    (\theta_1, \theta_1)_T & (\theta_2, \theta_1)_T & (\theta_3, \theta_1)_T \\
    (\theta_1, \theta_2)_T& (\theta_2, \theta_2)_T& (\theta_3, \theta_2)_T \\
    (\theta_1, \theta_3)_T& (\theta_2, \theta_3)_T& (\theta_3, \theta_3)_T\\
  \end{array}
\right)
\left(\begin{array}{c}c_1 \\c_2 \\c_3 \\\end{array}\right)=\left(\begin{array}{c}
-(v_0,\nabla\cdot\theta_1)_T+(v_b,\;\theta_1\cdot\bn)_{\partial T} \\
 -(v_0,\nabla\cdot\theta_2)_T+(v_b,\;\theta_2\cdot\bn)_{\partial T}\\
 -(v_0,\nabla\cdot\theta_3)_T+(v_b,\;\theta_3\cdot\bn)_{\partial T}
 \end{array}
 \right).
 \]
The inverse of the above coefficient matrix can be obtained
explicitly or numerically through a local matrix solver. Once
$\nabla_d\phi_i$ and $\nabla_d\psi_j$ are computed, the rest of
steps are the same as those for the classical finite element method.

\section{Error Estimates}

Denote by $Q_h u=\{Q_0 u,\;Q_bu\}$ the $L^2$ projection onto
$P_j(T^0)\times P_{j}(\partial T)$. In other words, on each element
$T$, the function $Q_0 u$ is defined as the $L^2$ projection of $u$
in $P_j(T)$ and on $\partial T$, $Q_b u$ is the $L^2$ projection in
$P_{j}(\partial T)$.

Consider the following Helmholtz problem:
\begin{eqnarray}
-\nabla\cdot(d\nabla u)-k^2u &=& f\quad
\mbox{in}\;\Omega,\label{pde1}\\
u&=&g\quad \mbox{on}\; \partial\Omega.\label{bc1}
\end{eqnarray}

For weak Galerkin finite element space $V_h$ defined in (\ref{vh}),
we define $V_h^0$ a subspace of $V_h$ with vanishing boundary values
on $\partial\Omega$; i.e.,
\[
V_h^0 :=\left\{ v=\{v_0, v_b\}\in V_h,
{v_b}|_{\partial T\cap \partial\Omega}=0, \ \forall T\in {\cal T}_h
\right\}.
\]
The weak Galerkin finite element method for the Helmholtz problem
(\ref{pde1})-(\ref{bc1}) is stated as follows.

\begin{algorithm}
A numerical approximation for (\ref{pde1}) and (\ref{bc1}) can be
obtained by seeking $u_h=\{u_0,u_b\}\in V_h$ satisfying
$u_b= Q_b g$ on $\partial \Omega$ and the following equation:
\begin{equation}\label{WG-fem}
(d\nabla_du_h,\;\nabla_dv)-k^2(u_0,\;v_0)=(f,\;v_0), \quad\forall\
v=\{v_0, v_b\}\in V_h^0,
\end{equation}
\end{algorithm}

For a sufficiently small mesh size $h$, the following error estimate
holds true. A proof of this error estimate can be found in
\cite{wy}.

\medskip

\begin{theorem}\label{H1error-estimate}
Let $u\in H^{1}(\Omega)$ be the solution
(\ref{pde1}) and (\ref{bc1}), and $u_h$ be a weak Galerkin
approximation of $u$ arising from (\ref{WG-fem}). Assume that
the exact solution $u$ is sufficiently smooth such that $u\in
H^{m+1}(\Omega)$ with $0\le m \le k+1$. Then, there exists a
constant $C$ such that
\begin{eqnarray}
\|\nabla_d(u_h-Q_hu)\|& \le& C(h^{m} \|u\|_{m+1}
+h^{1+s}\|f-Q_0f\|)\label{H1error}\\
\|u_h-Q_hu\|  &\le& C(h^{m+s} \|u\|_{m+1}
+h^{1+s}\|f-Q_0f\|),\label{L2error}
\end{eqnarray}
provided that the dual problem of (\ref{pde1})-(\ref{bc1}) has the
$H^{m+s}(\Omega)$ regularity with $s\in (0,1]$.
\end{theorem}

\section{Numerical Experiments}
In this section, we examine the WG method by testing its accuracy,
convergence, and robustness for solving two dimensional Helmholtz
equations. The pollution effect due to large wave numbers will be
particularly investigated and tested numerically. For convergence
tests, both piecewise constant and piecewise linear finite elements
will be considered. To demonstrate WG's robustness, the Helmholtz
equation in both homogeneous and inhomogeneous media will be solved
on convex and non-convex computational domains. For simplicity, a
structured triangular mesh is employed in all cases. The mesh
generation and all computations are conducted in the MATLAB
environment.

Two types of relative errors are measured in our numerical
experiments. The first one is the relative $L^2$ error and the
second is the relative $H^1$ error. The relative $L^2$ error is
defined by
$$ \frac{\| u_h - Q_h u \|}{\| Q_h u \|},$$
where $u_h$ is the WG approximation to the solution $u$, and $Q_h u
= \{ Q_0 u, Q_b u \}$ is the $L^2$ projection onto $P_j(T^0) \times
P_j(\partial T)$. In other words, on each element $T$, the function
$Q_0 u$ is defined as the $L^2$ projection of $u$ in $P_j(T)$ and on
$\partial T$, $Q_b u$ is the $L^2$ projection in $P_{j}(\partial
T)$. The relative $H^1$ error is defined in terms of the discrete
gradient
$$ \frac{\| \nabla_d (u_h - Q_h u) \|}{\| \nabla_d Q_h u \|}.$$
However, in the present study, the $H^1$-semi-norm will be calculated
as
$$
\3bar u_h-Q_hu \3bar^2=h_e^{-1}\langle
u_0-u_b-(Q_0u-Q_bu),u_0-u_b-(Q_0u-Q_bu)\rangle_e
$$
for the lowest order finite element (i.e., piecewise constants). For
piecewise linear elements, we use the original definition of
$\nabla_d$ to compute the $H^1$-semi-norm $\|\nabla_d (u_h - Q_h
u)\|$.

\begin{figure}[!tb]
\centering
\begin{tabular}{cc}
  \resizebox{2.25in}{2.25in}{\includegraphics{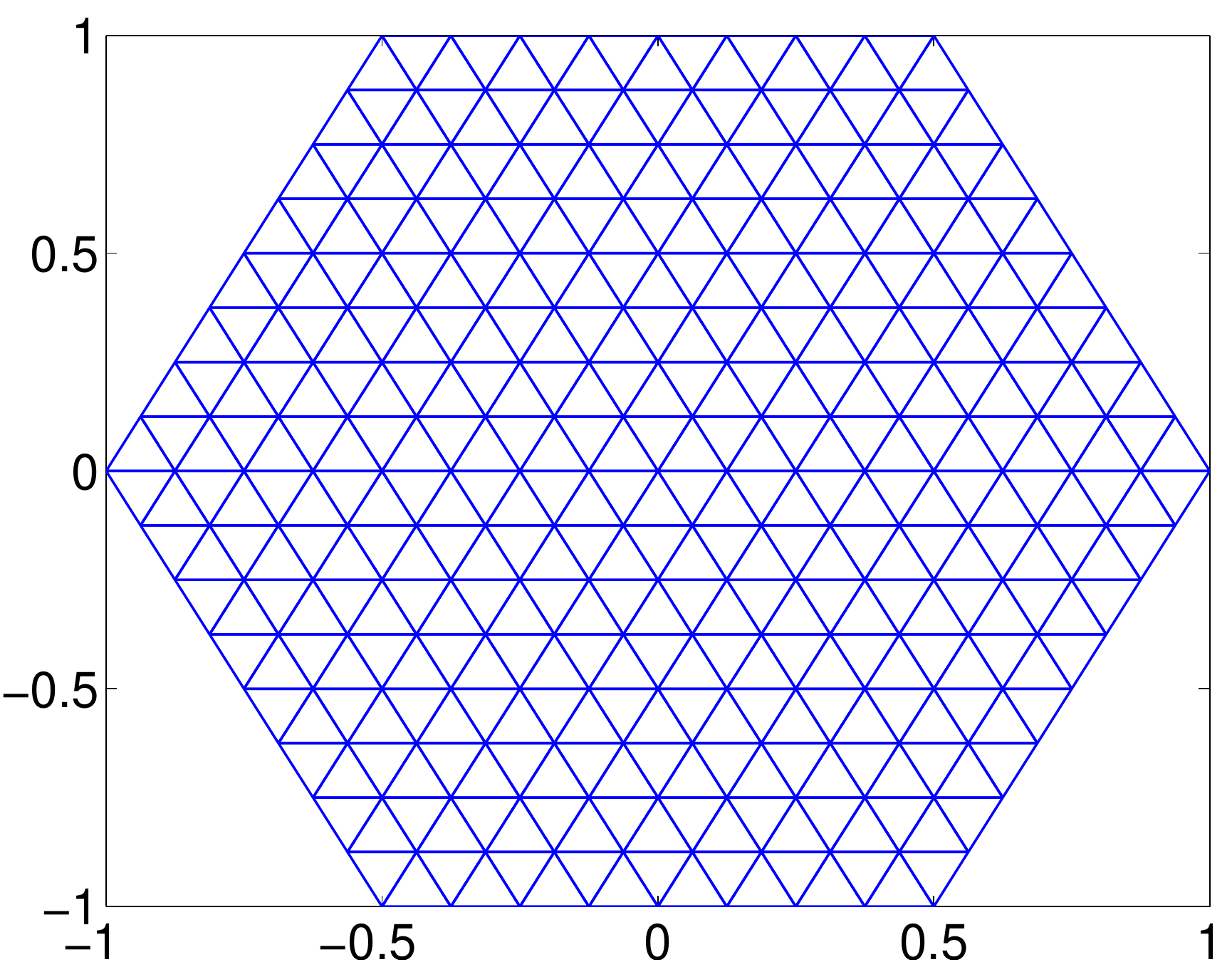}} \quad
  \resizebox{2.25in}{2.25in}{\includegraphics{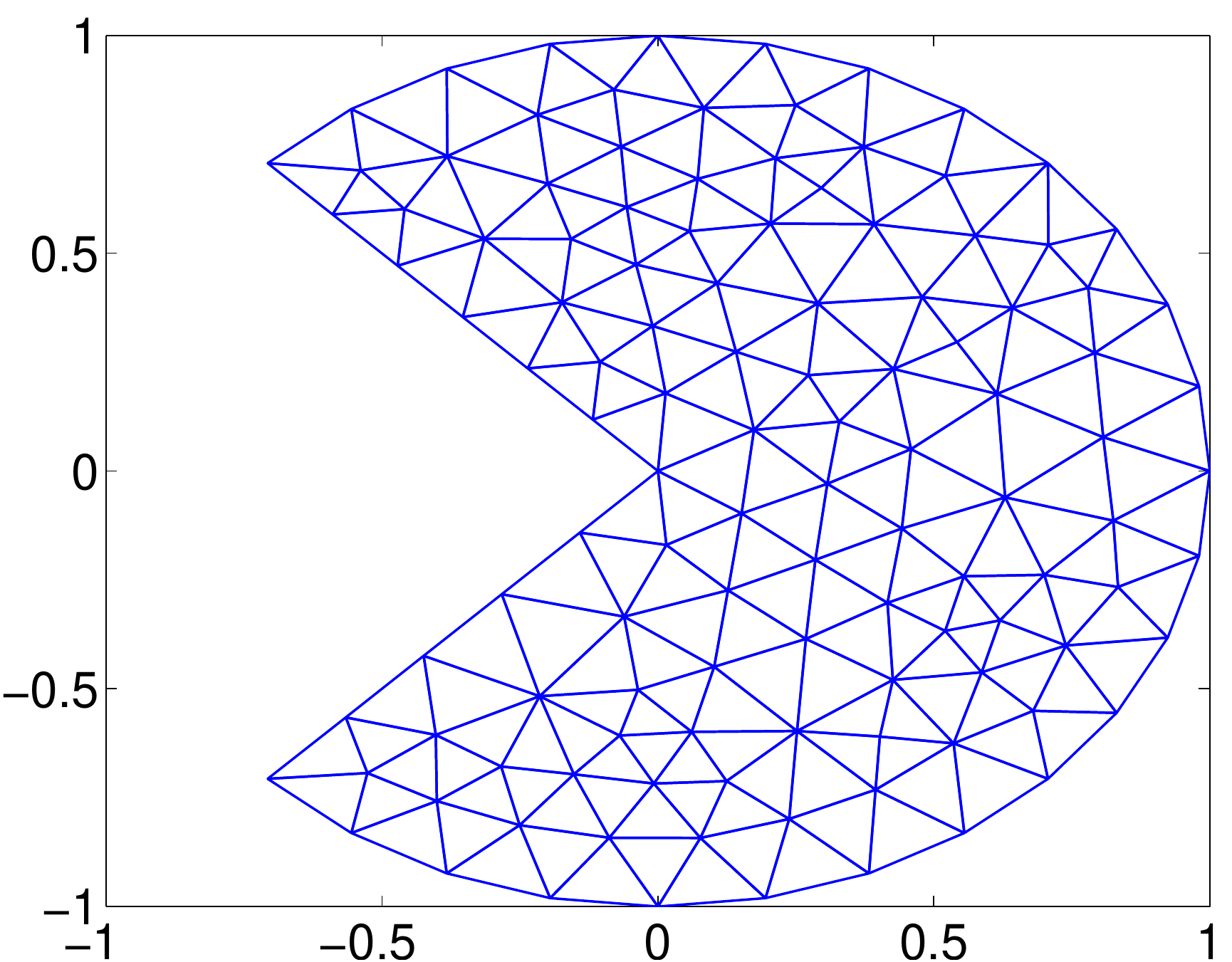}}
\end{tabular}
\caption{Geometry of testing domains and sample meshes. Left: a
convex hexagon domain; Right: a non-convex imperfect circular
domain.} \label{fig.domain}
\end{figure}

\subsection{A convex Helmholtz problem}\label{convex}
We first consider a homogeneous Helmholtz equation defined on a
convex hexagon domain, which has been studied in \cite{fw}. The
domain $\Omega$ is the unit regular hexagon domain centered at the
origin $(0,0)$, see Fig. \ref{fig.domain} (left). Here we set $d=1$
and $f=\sin(kr)/r$ in (\ref{pde}), where $r=\sqrt{x^2+y^2}$. The
boundary data $g$ in the Robin boundary condition (\ref{bc}) is
chosen so that the exact solution is given by
\begin{eqnarray}
u=\frac{\cos(kr)}{k}-\frac{\cos k+i\sin k}{k(J_0(k)+iJ_1(k))}J_0(kr)
\end{eqnarray}
where $J_{\xi}(z)$ are Bessel functions of the first kind. Let
$\T_h$ denote the regular triangulation that consists of $6N^2$
triangles of size $h=1/N$, as shown in Fig. \ref{fig.domain} (left)
for $T_{\frac18}$.

\begin{table}[!t]
\caption{Convergence of piecewise constant WG for the Helmholtz equation on a
convex domain with wave number $k=1$.} \label{table.Ex1k1}
\begin{center}
\begin{tabular}{|l|l|l|l|l|}
\hline
& \multicolumn{2}{c|}{relative $H^1$} & \multicolumn{2}{c|}{relative $L^2$} \\
\cline{2-3} \cline{4-5}
$h$ & error & order & error & order \\
\hline
  5.00e-01 &  2.49e-02 &      &  4.17e-03 &\\
  2.50e-01 &  1.11e-02 & 1.16 &  1.05e-03 & 1.99 \\
  1.25e-01 &  5.38e-03 & 1.05 &  2.63e-04 & 2.00 \\
  6.25e-02 &  2.67e-03 & 1.01 &  6.58e-05 & 2.00 \\
  3.13e-02 &  1.33e-03 & 1.00 &  1.64e-05 & 2.00 \\
  1.56e-02 &  6.65e-04 & 1.00 &  4.11e-06 & 2.00 \\
\hline
\end{tabular}
\end{center}
\end{table}

Table \ref{table.Ex1k1} illustrates the performance of the WG method
with piecewise constant elements for the Helmholtz equation with
wave number $k=1$. Uniform triangular partitions were used in the
computation through successive mesh refinements. The relative errors
in $L^2$ norm and $H^1$ semi-norm can be seen in Table
\ref{table.Ex1k1}. The Table also includes numerical estimates for
the rate of convergence in each metric. It can be seen that the
order of convergence in the relative $H^1$ semi-norm and relative
$L^2$ norm are, respectively, one and two for piecewise constant
elements.

\begin{table}[!t]
\caption{ Convergence of piecewise linear WG for the Helmholtz equation on a
convex domain with wave number $k=5$.} \label{table.Ex1k5}
\begin{center}
\begin{tabular}{|l|l|l|l|l|}
\hline
& \multicolumn{2}{c|}{relative $H^1$} & \multicolumn{2}{c|}{relative $L^2$} \\
\cline{2-3} \cline{4-5}
$h$ & error & order & error & order \\
\hline
  2.50e-01 &  9.48e-03 &      & 2.58e-04 & \\
  1.25e-01 &  2.31e-03 & 2.04 & 3.46e-05 & 2.90 \\
  6.25e-02 &  5.74e-04 & 2.01 & 4.47e-06 & 2.95 \\
  3.13e-02 &  1.43e-04 & 2.00 & 5.64e-07 & 2.99 \\
  1.56e-02 &  3.58e-05 & 2.00 & 7.06e-08 & 3.00 \\
  7.81e-03 &  8.96e-06 & 2.00 & 8.79e-09 & 3.01 \\
\hline
\end{tabular}
\end{center}
\end{table}

High order of convergence can be achieved by using corresponding
high order finite elements in the present WG framework. To
demonstrate this phenomena, we consider the same Helmholtz problem
with a slightly larger wave number $k=5$. The WG with piecewise
linear functions was employed in the numerical approximation. The
computational results are reported in Table \ref{table.Ex1k5}. It is
clear that the numerical experiment validates the theoretical
estimates. More precisely, the rates of convergence in the relative
$H^1$ semi-norm and relative $L^2$ norm are given by two and three,
respectively.

\subsection{A non-convex Helmholtz problem}

We next explore the use of the WG method for solving a Helmholtz
problem defined on a non-convex domain, see Fig. \ref{fig.domain}
(right). The medium is still assumed to be homogeneous, i.e., $d=1$
in (\ref{pde}). We are particularly interested in the performance of
the WG method for dealing with the possible field singularity at the
origin. For simplicity, only the piecewise constant $RT_0$ elements
are tested for the present problem. Following \cite{Monk11}, we take
$f=0$ in (\ref{pde}) and the boundary condition is simply taken as a
Dirichlet one: $u=g$ on $\partial \Omega$. Here $g$ is prescribed
according to the exact solution \cite{Monk11}
\begin{equation}\label{solution2}
u= J_{\xi}( k \sqrt{x^2+y^2}) \cos (\xi \arctan( y/x)).
\end{equation}

\begin{figure}[!tb]
\centering
\begin{tabular}{cc}
  \resizebox{2.55in}{2.1in}{\includegraphics{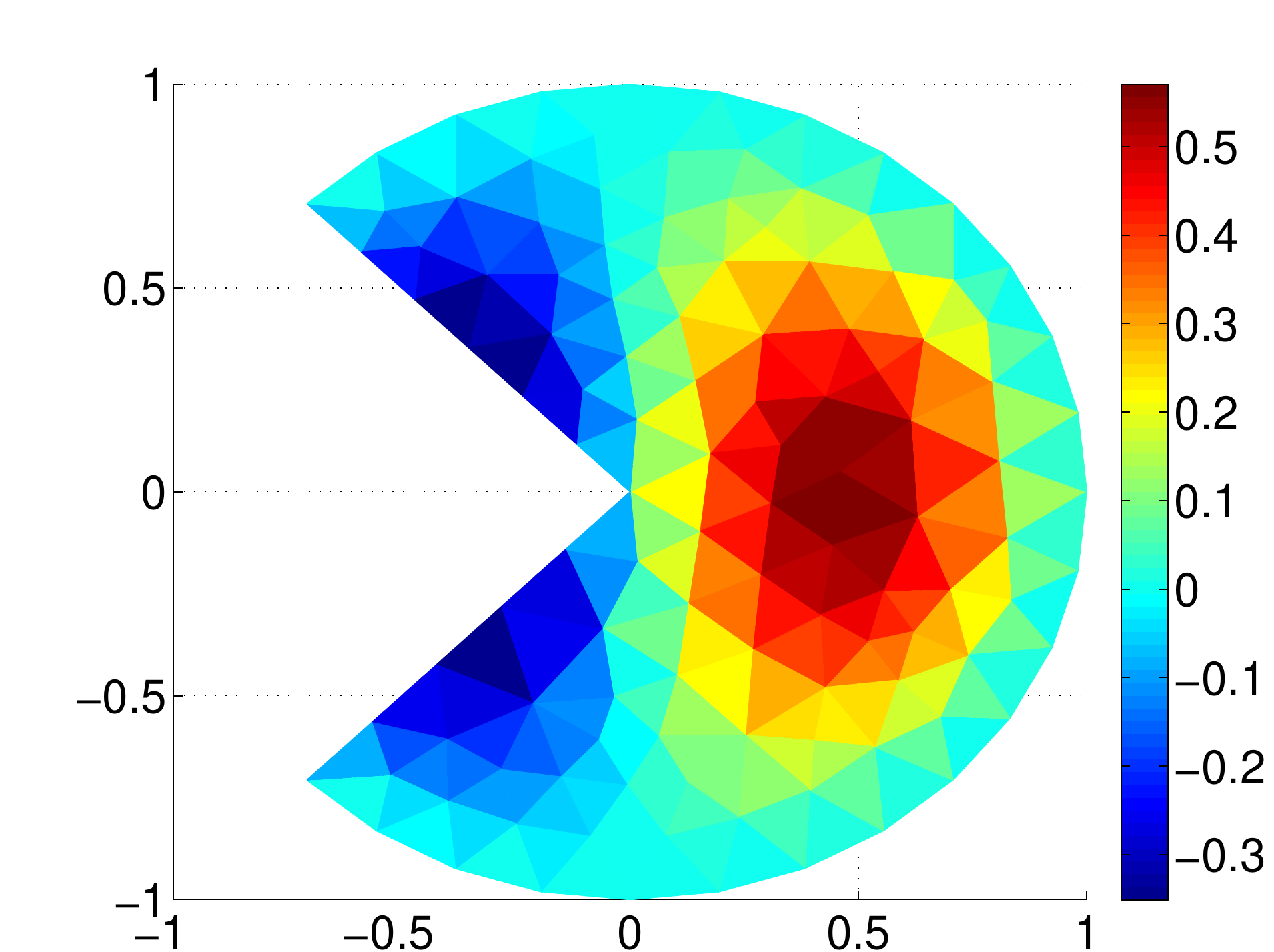}} \quad
  \resizebox{2.45in}{2.1in}{\includegraphics{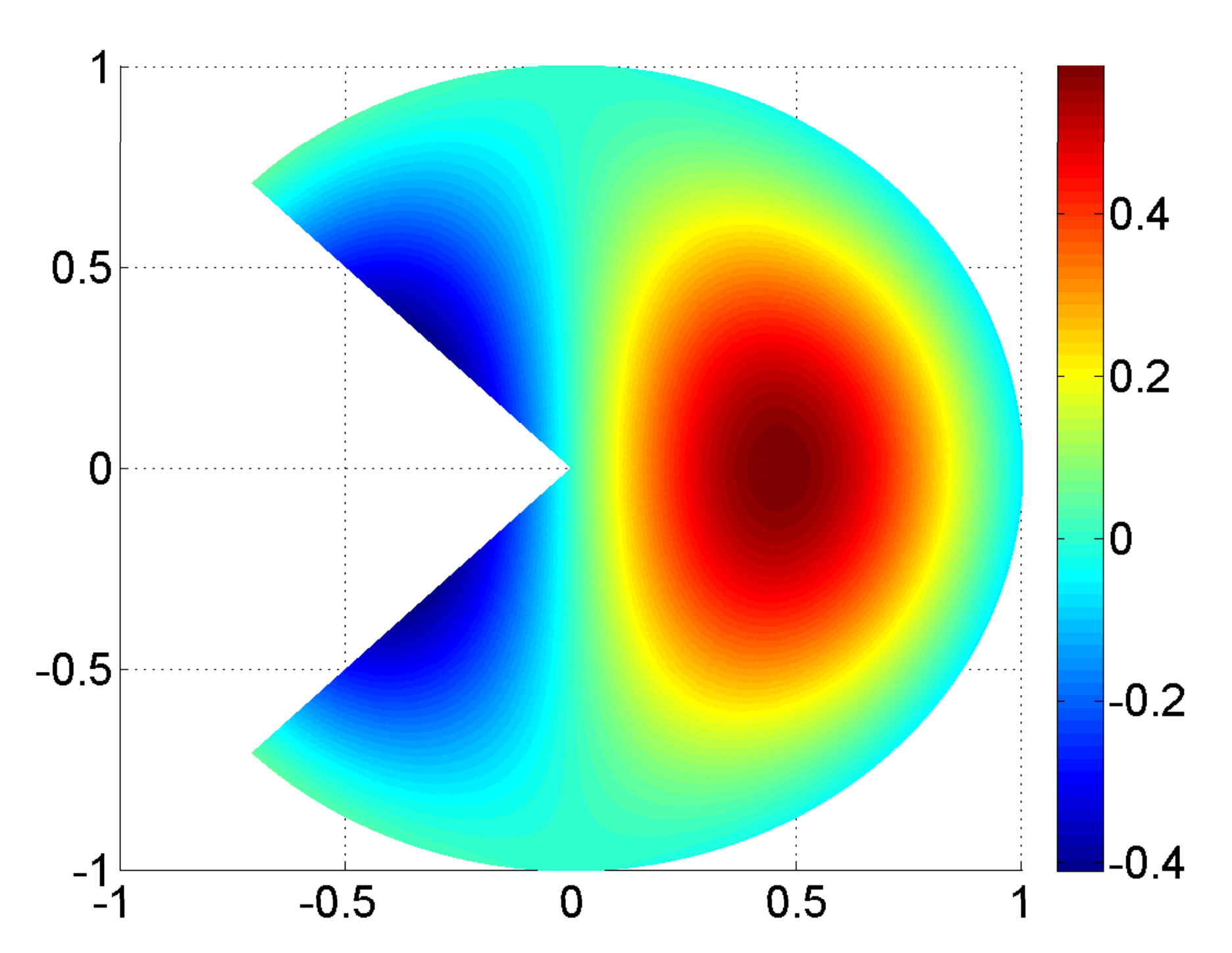}}
\end{tabular}
\caption{WG solutions for the non-convex Helmholtz problem with
$k=4$ and $\xi=1$. Left: Mesh level $1$; Right: Mesh level $6$.}
\label{fig.Ex2xi1}
\end{figure}

\begin{figure}[!tb]
\centering
\begin{tabular}{cc}
  \resizebox{2.55in}{2.1in}{\includegraphics{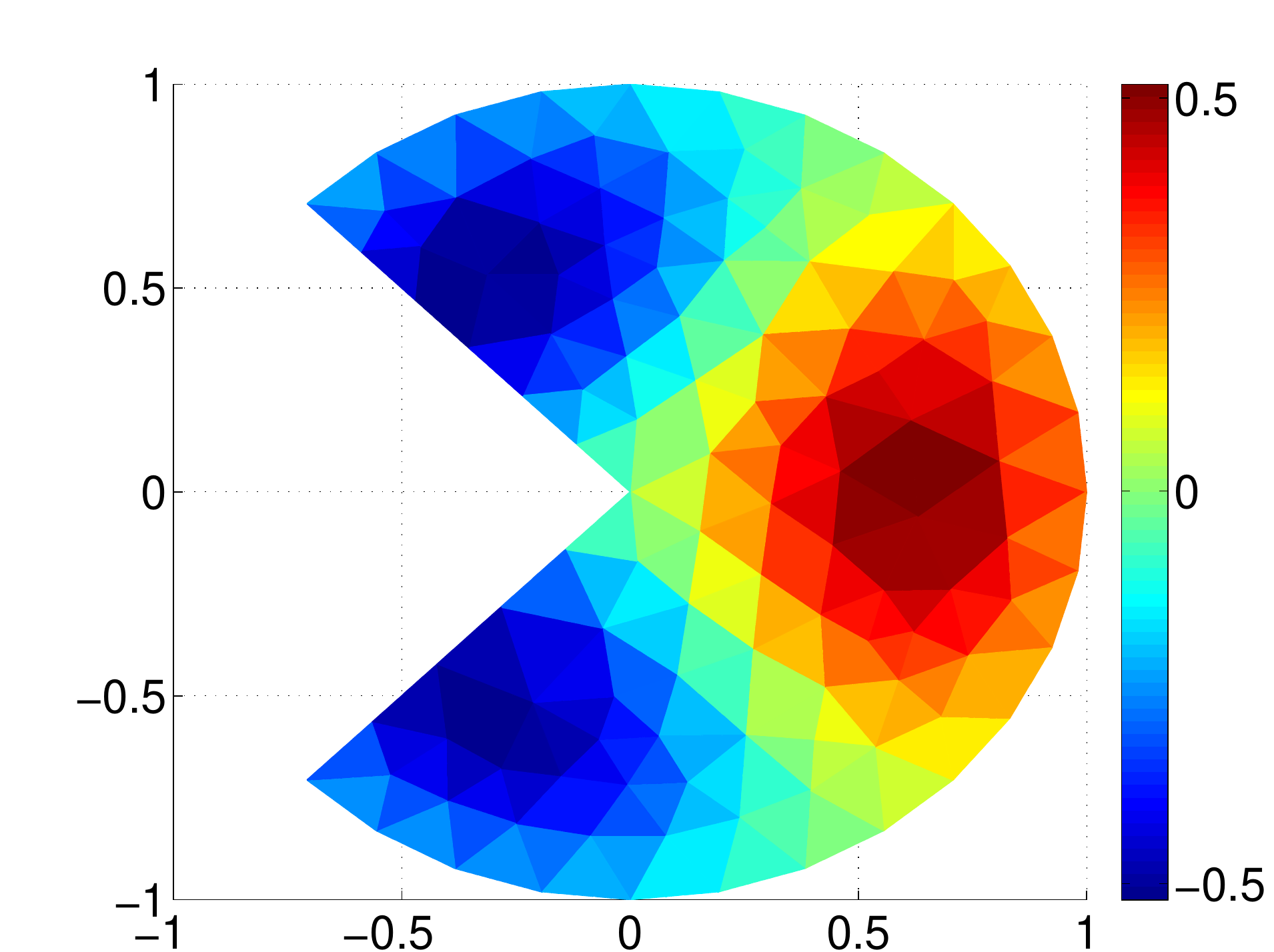}} \quad
  \resizebox{2.45in}{2.1in}{\includegraphics{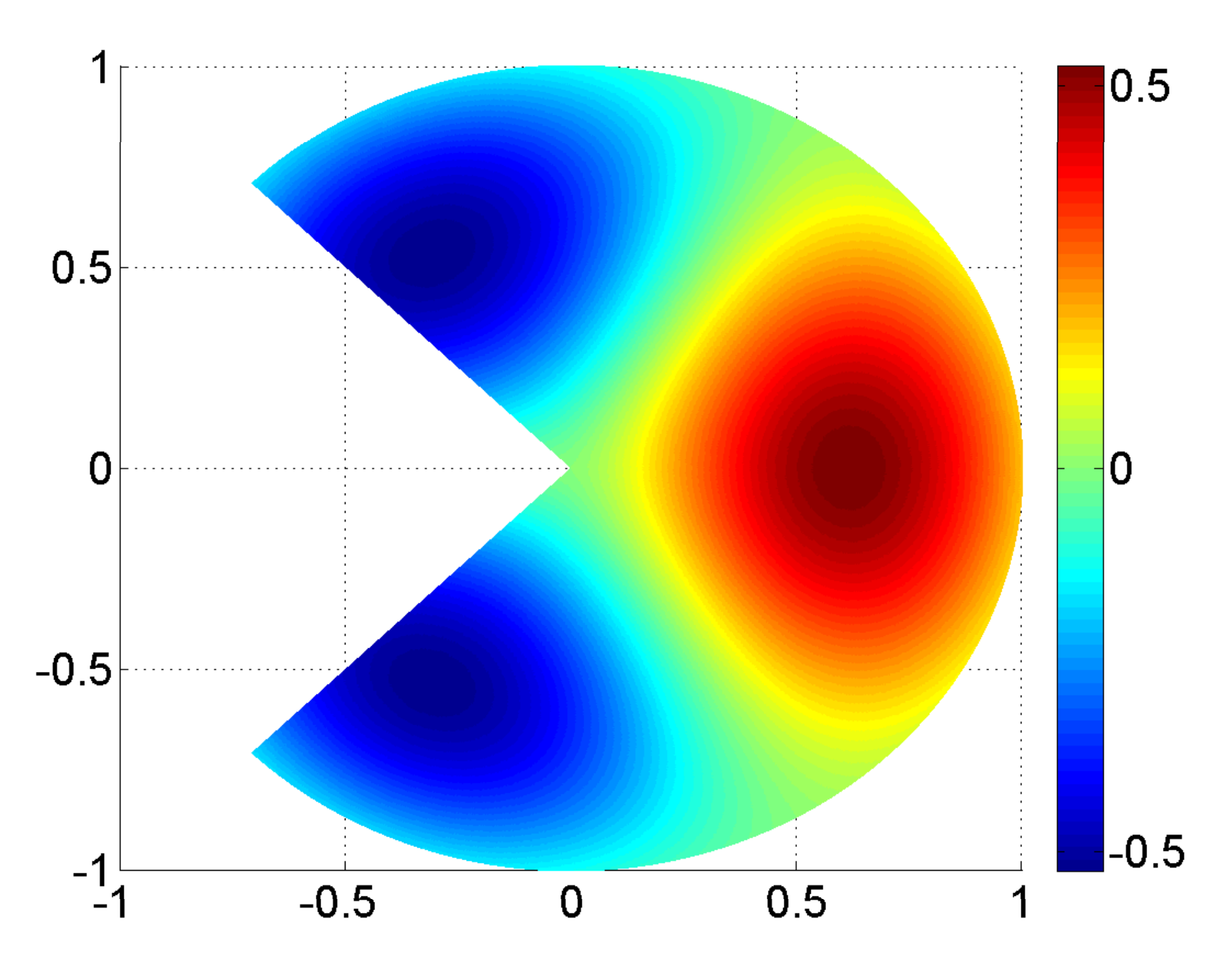}}
\end{tabular}
\caption{WG solutions for the non-convex Helmholtz problem with
$k=4$ and $\xi=3/2$. Left: Mesh level $1$; Right: Mesh level $6$.}
\label{fig.Ex2xi32}
\end{figure}

\begin{figure}[!tb]
\centering
\begin{tabular}{cc}
  \resizebox{2.55in}{2.1in}{\includegraphics{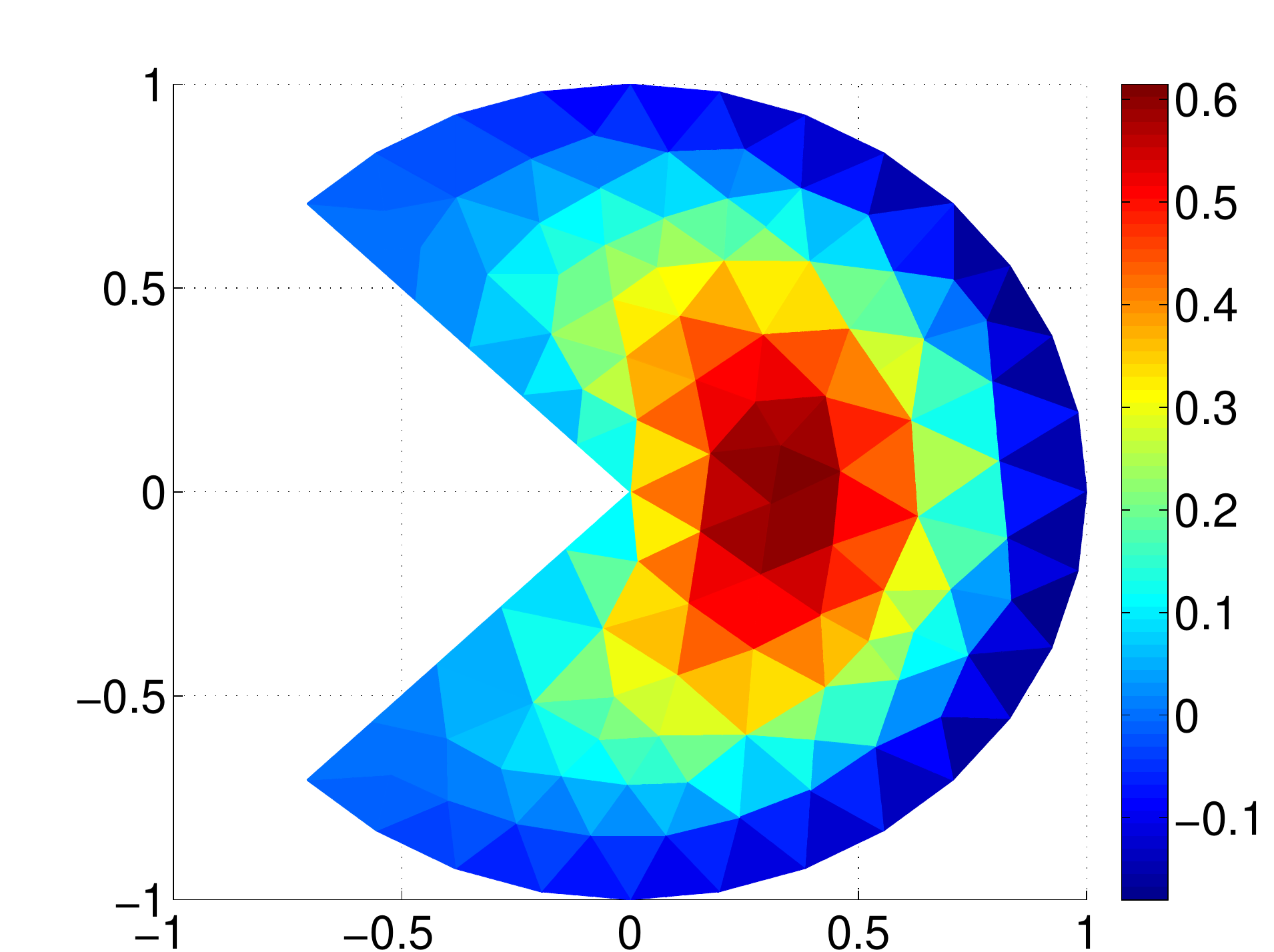}}
  \resizebox{2.45in}{2.1in}{\includegraphics{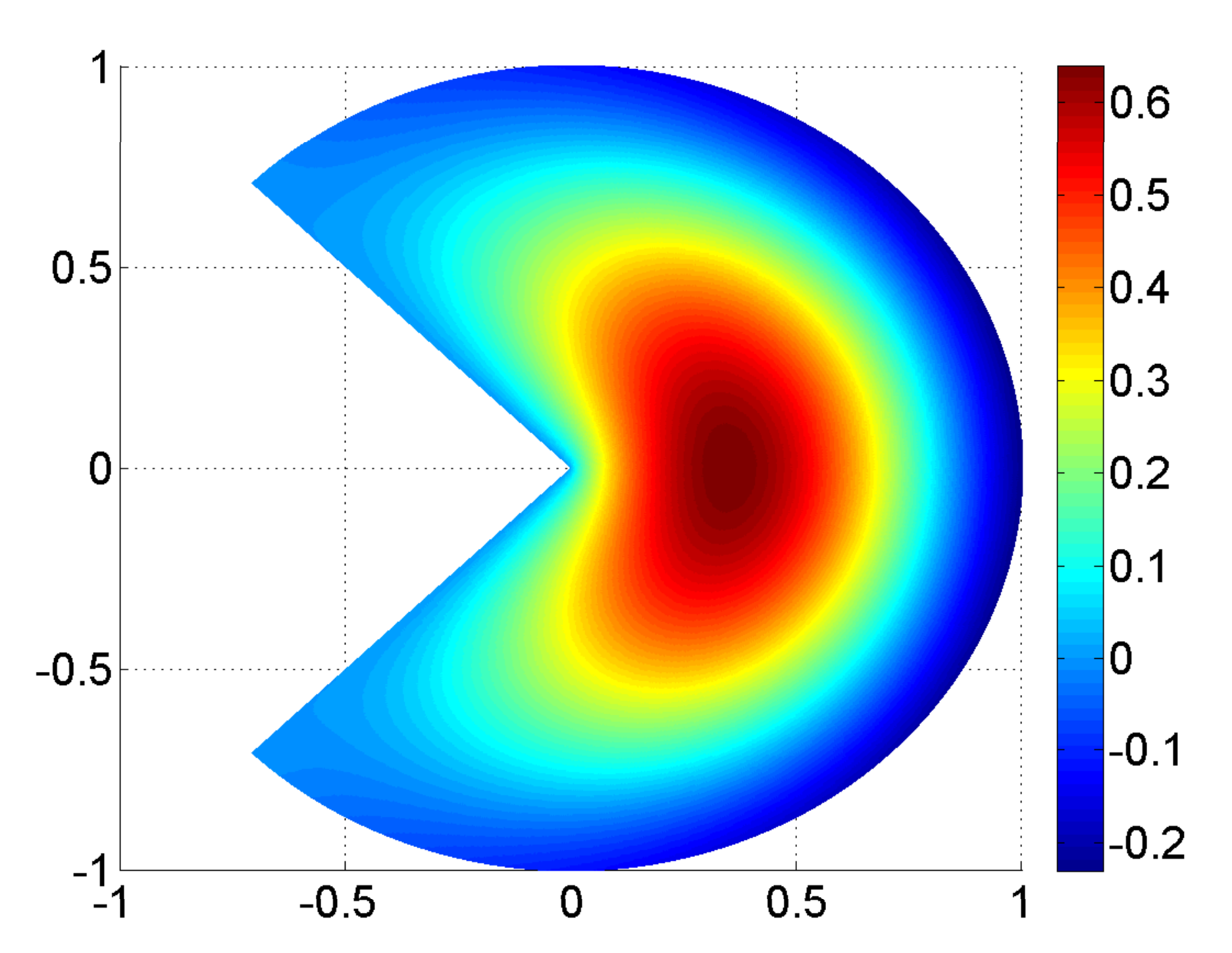}}
\end{tabular}
\caption{WG solutions for the non-convex Helmholtz problem with
$k=4$ and $\xi=2/3$. Left: Mesh level $1$; Right: Mesh level $6$.}
\label{fig.Ex2xi23}
\end{figure}

In the present study, the wave number was chosen as $k=4$ and three
values for the parameter $\xi$ are considered; i.e., $\xi=1$,
$\xi=3/2$ and $\xi=2/3$. The same triangular mesh is used in the WG
method for all three cases. In particular, an initial mesh is first
generated by using MATLAB with default settings, see Fig.
\ref{fig.domain} (right). Next, the mesh is refined uniformly for
five times. The WG solutions on mesh level $1$ and mesh level $6$
are shown in Fig. \ref{fig.Ex2xi1}, Fig. \ref{fig.Ex2xi32}, and Fig.
\ref{fig.Ex2xi23}, respectively, for $\xi=1$, $\xi=3/2$, and
$\xi=2/3$. Since the numerical errors are quite small for the WG
approximation corresponding to mesh level $6$, the field modes
generated by the densest mesh are visually indistinguishable from
the analytical ones. In other words, the results shown in the right
charts of Fig. \ref{fig.Ex2xi1}, Fig. \ref{fig.Ex2xi32}, and Fig.
\ref{fig.Ex2xi23} can be regarded as analytical results. It can be
seen that in all three cases, the WG solutions already agree with
the analytical ones at the coarsest level. Moreover, based on the
coarsest mesh, the constant function values can be clearly seen in
each triangle, due to the use of piecewise constant $RT_0$ elements.
Nevertheless, after the initial mesh is refined for five times, the
numerical plots shown in the right charts are very smooth. A perfect
symmetry with respect to the $x$-axis is clearly seen.

\begin{table}[!tb]
\caption{Numerical convergence test for the non-convex Helmholtz
problem with $k=4$ and $\xi=1$. } \label{table.Ex2xi1}
\begin{center}
\begin{tabular}{|l|l|l|l|l|}
\hline
& \multicolumn{2}{c|}{relative $H^1$} & \multicolumn{2}{c|}{relative $L^2$} \\
\cline{2-3} \cline{4-5}
$h$ & error & order & error & order \\
\hline
  2.44e-01 & 5.64e-02 &      & 1.37e-02 & \\
  1.22e-01 & 2.83e-02 & 1.00 & 3.56e-03 & 1.95 \\
  6.10e-02 & 1.42e-02 & 0.99 & 8.98e-04 & 1.99 \\
  3.05e-02 & 7.14e-03 & 1.00 & 2.25e-04 & 2.00 \\
  1.53e-02 & 3.57e-03 & 1.00 & 5.63e-05 & 2.00 \\
  7.63e-03 & 1.79e-03 & 1.00 & 1.41e-05 & 2.00 \\
\hline
\end{tabular}
\end{center}
\end{table}

\begin{table}[!tb]
\caption{Numerical convergence test for the non-convex Helmholtz
problem with $k=4$ and $\xi=3/2$. } \label{table.Ex2xi32}
\begin{center}
\begin{tabular}{|l|l|l|l|l|}
\hline
& \multicolumn{2}{c|}{relative $H^1$} & \multicolumn{2}{c|}{relative $L^2$} \\
\cline{2-3} \cline{4-5}
$h$ & error & order & error & order \\
\hline
  2.44e-01   & 5.56e-02 & & 1.12e-2& \\
  1.22e-01   & 2.81e-02 & 0.98 & 3.02e-03 & 1.89 \\
  6.10e-02   & 1.42e-02 & 0.99 & 8.06e-04 & 1.91 \\
  3.05e-02   & 7.14e-03 & 0.99 & 2.12e-04 & 1.92 \\
  1.53e-02   & 3.58e-03 & 1.00 & 5.54e-05 & 1.94 \\
  7.63e-03   & 1.79e-03 & 1.00 & 1.44e-05 & 1.95 \\
\hline
\end{tabular}
\end{center}
\end{table}

\begin{table}[!tb]
\caption{Numerical convergence test for the non-convex Helmholtz
problem with $k=4$ and $\xi=2/3$. } \label{table.Ex2xi23}
\begin{center}
\begin{tabular}{|l|l|l|l|l|}
\hline
& \multicolumn{2}{c|}{relative $H^1$} & \multicolumn{2}{c|}{relative $L^2$} \\
\cline{2-3} \cline{4-5}
$h$ & error & order & error & order \\
\hline
   2.44e-01   & 1.07e-01 &       &5.24e-02 & \\
   1.22e-01   & 5.74e-02 & 0.90  &2.18e-02 & 1.27 \\
   6.10e-02   & 3.23e-02 & 0.83  &9.01e-03 & 1.27 \\
   3.05e-02   & 1.89e-02 & 0.77  &3.68e-03 & 1.29 \\
   1.53e-02   & 1.14e-02 & 0.73  &1.49e-03 & 1.31 \\
   7.63e-03   & 6.99e-03 & 0.71  &5.96e-04 & 1.32 \\
\hline
\end{tabular}
\end{center}
\end{table}

\begin{figure}[!tb]
\centering
\epsfig{file=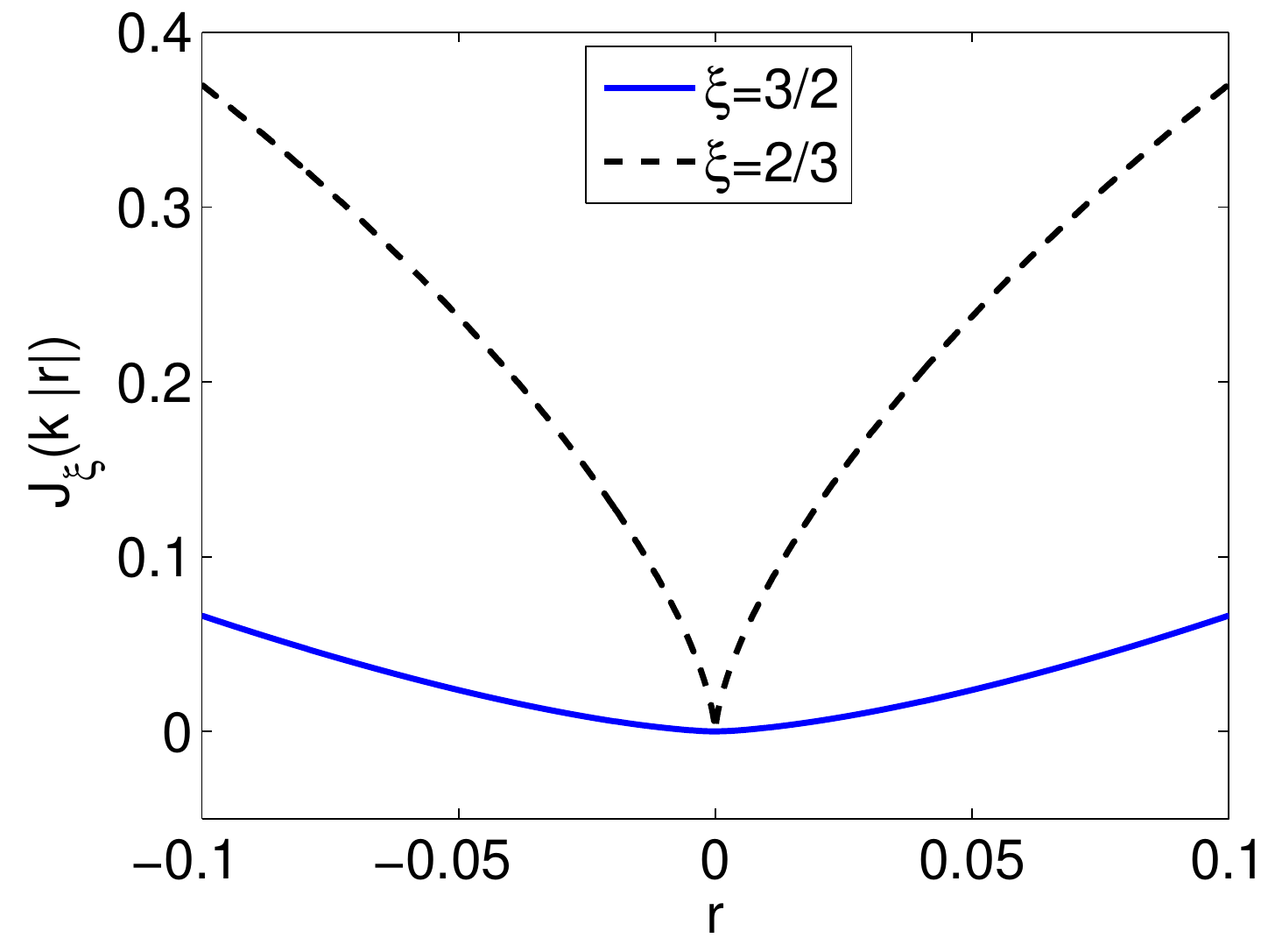,width=2.5in}
\caption{The Bessel function of the first kind
$J_{\xi}( k |r|)$ across the origin. }
\label{fig.origin}
\end{figure}

We next investigate the numerical convergence rates for WG. The
numerical errors of the WG solutions for $\xi=1$, $\xi=3/2$ and
$\xi=2/3$ are listed, respectively, in Table \ref{table.Ex2xi1},
Table \ref{table.Ex2xi32}, and Table \ref{table.Ex2xi23}. It can be
seen that for $\xi=1$ and $\xi=3/2$, the numerical convergence rates
in the relative $H^1$ and $L^2$ errors remain to be first and second
order, while the convergence orders degrade for the non-smooth case
$\xi=2/3$. Mathematically, for both  $\xi=3/2$ and $\xi=2/3$, the
exact solutions (\ref{solution2}) are known to be non-smooth across
the negative $x$-axis if the domain was chosen to be the entire
circle. However, the present domain excludes the negative $x$-axis.
Thus, the source term $f$ of the Helmholtz equation (\ref{pde}) can
be simply defined as zero throughout $\Omega$. Nevertheless, there
still exists some singularities at the origin $(0,0)$. In
particular, it is remarked in \cite{Monk11} that the singularity
lies in the derivatives of the exact solution at $(0,0)$. Due to
such singularities, the convergence rates of high order DG methods
are also reduced for $\xi=3/2$ and $\xi=2/3$ \cite{Monk11}. In the
present study, we further note that there exists a subtle difference
between two cases $\xi=3/2$ and $\xi=2/3$ at the origin. To see
this, we neglect the second $\cos(\cdot)$ term in the exact solution
(\ref{solution2}) and plot the Bessel function of the first kind
$J_{\xi}( k |r|)$ along the radial direction $r$, see Fig.
\ref{fig.origin}. It is observed that the Bessel function of the
first kind is non-smooth for the case $\xi=2/3$, while it looks
smooth across the origin for the case $\xi=3/2$. Thus, it seems that
the first derivative of $J_{3/2}( k |r|)$ is still continuous along
the radial direction. This perhaps explains why the present WG
method does not experience any order reduction for the case
$\xi=3/2$. In  \cite{Monk11}, locally refined meshes were employed
to resolve the singularity at the origin so that the convergence
rate for the case $\xi=2/3$ can be improved. We note that local
refinements can also be adopted in the WG method for a better
convergence rate. A study of WG with grid local refinement is left
to interested parties for future research.

\subsection{A Helmholtz problem with inhomogeneous media}

\begin{figure}[!tb]
\centering
\resizebox{2.5in}{2.15in}{\includegraphics{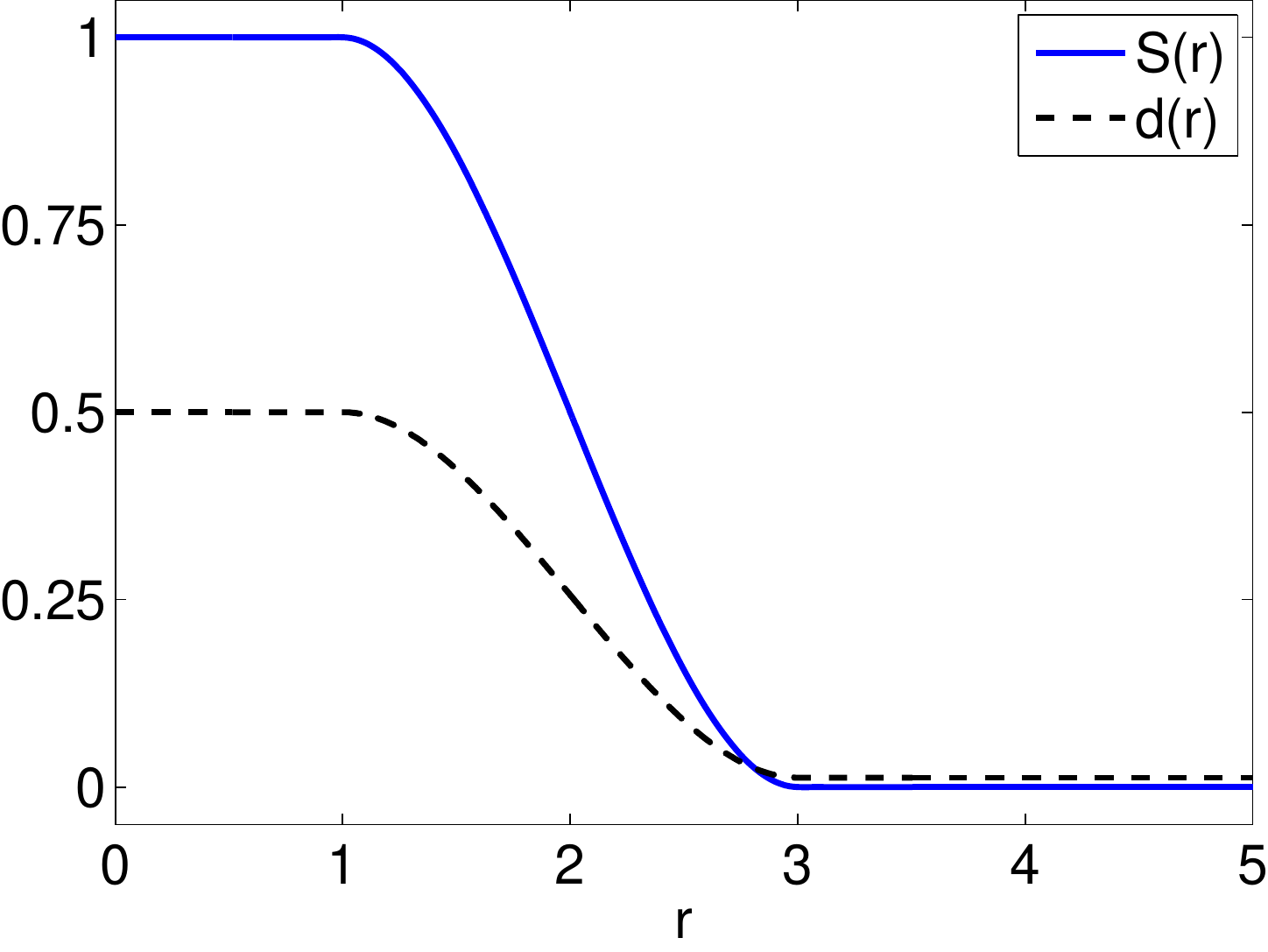}}
\caption{An example plot of smooth dielectric profile
$d(r)$ and $S(r)$ with $a=1$, $b=3$ and $R=5$. The dielectric
coefficients of  protein and water are used,
i.e., $\epsilon_1=2$ and $\epsilon_2=80$. }
\label{fig.eps}
\end{figure}

We consider a Helmholtz problem with inhomogeneous media defined on
a circular domain with radius $R$. Note that the spatial function
$d(x,y)$ in the Helmholtz equation (\ref{pde}) represents the
dielectric properties of the underlying media. In particular, we
have $d=\frac{1}{\epsilon}$ in the electromagnetic applications
\cite{Zhao10}, where $\epsilon$ is the electric permittivity. In the
present study, we construct a smooth varying dielectric profile:
\begin{equation}\label{dr}
d(r)= \frac{1}{\epsilon_1}S(r) + \frac{1}{\epsilon_2}(1-S(r)),
\end{equation}
where $r=\sqrt{x^2+y^2}$, $\epsilon_1$ and $\epsilon_2$ are
dielectric constants, and
\begin{equation}
S(r)=
\begin{cases}
1 & \text{if $r<a$}, \\
-2\left( \frac{b-r}{b-a} \right)^3 +
 3\left( \frac{b-r}{b-a} \right)^2 & \text{if $a \le r \le b$}, \\
0 & \text{if $r>b$},
\end{cases}
\end{equation}
with $a<b<R$. An example plot of $d(r)$ and $S(r)$ is shown in Fig.
\ref{fig.eps}. In classical electromagnetic simulations, $\epsilon$
is usually taken as a piecewise constant, so that  some
sophisticated numerical treatments have to be conducted near the
material interfaces to secure the overall accuracy \cite{Zhao10}.
Such a procedure can be bypassed if one considers a smeared
dielectric profile, such as (\ref{dr}). We note that under the limit
$b \to a$, a piecewise constant profile is recovered in (\ref{dr}).
In general, the smeared profile (\ref{dr}) might be generated via
numerical filtering, such as the so-called $\epsilon$-smoothing
technique \cite{Shao03} in computational electromagnetics. On the
other hand, we note that the dielectric profile might be defined to
be smooth in certain applications. For example, in studying the
solute-solvent interactions of electrostatic analysis, some
mathematical models \cite{Chen10,Zhao11} have been proposed to treat
the boundary between the protein and its surrounding aqueous
environment to be a smoothly varying one. In fact, the definition of
(\ref{dr}) is inspired by a similar model in that field
\cite{Chen10}.

In the present study, we choose the source of the Helmholtz equation
(\ref{pde}) to be
\begin{equation}
f(r)=k^2 [d(r)-1] J_0 (k r) + k d'(r) J_1(kr),
\end{equation}
where
\begin{equation}
d'(r)=\left( \frac{1}{\epsilon_1}-\frac{1}{\epsilon_2} \right) S'(r)
\end{equation}
and
\begin{equation}
S'(r)=
\begin{cases}
0 & \text{if $r<a$}, \\
6\left( \frac{b-r}{b-a} \right)^2 -
6\left( \frac{b-r}{b-a} \right) & \text{if $a \le r \le b$}, \\
0 & \text{if $b<r$},
\end{cases}
\end{equation}
For simplicity, a Dirichlet boundary condition is imposed at $r=R$
with $u=g$. Here $g$ is prescribed according to the exact solution
\begin{equation}
u=J_0(kr).
\end{equation}

Our numerical investigation assumes the value of $a=1$, $b=3$ and
$R=5$. The wave number is set to be $k=2$. The dielectric
coefficients are chosen as $\epsilon_1=2$ and $\epsilon_2=80$, which
represents the dielectric constant of protein and water
\cite{Chen10,Zhao11}, respectively. The WG method with piecewise
constant finite element functions is employed to solve the present
problem with inhomogeneous media in Cartesian coordinate. Table
\ref{table.Ex3} illustrates the computational errors and some
numerical rate of convergence. It can be seen that the numerical
convergence in the relative $L^2$ error is not uniform, while the
relative $H^1$ error still converges uniformly in first order. This
phenomena might be related to the non-uniformity and smallness of
the media in part of the computational domain. Nevertheless, the
averaged convergence rate in the relative $L^2$ norm is about
$2.12$. Overall, we are confident that the WG method is accurate and
robust in solving the Helmholtz equations with inhomogeneous media.

\begin{table}[!h]
\caption{Numerical convergence test of the Helmholtz equation with
inhomogeneous media. } \label{table.Ex3}
\begin{center}
\begin{tabular}{|l|l|l|l|l|}
\hline
& \multicolumn{2}{c|}{relative $H^1$} & \multicolumn{2}{c|}{relative $L^2$} \\
\cline{2-3} \cline{4-5}
$h$ & error & order & error & order \\
\hline
   1.51e-00   & 2.20e-01 &       & 1.04e-00 &  \\
   7.54e-01   & 1.24e-01 & 0.83  & 1.20e-01 & 3.11 \\
   3.77e-01   & 6.24e-02 & 0.99  & 1.81e-02 & 2.73 \\
   1.88e-01   & 3.13e-02 & 1.00  & 5.71e-03 & 1.67 \\
   9.42e-02   & 1.56e-02 & 1.00  & 2.14e-03 & 1.42 \\
   4.71e-02   & 7.82e-03 & 1.00  & 5.11e-04 & 2.06 \\
\hline
\end{tabular}
\end{center}
\end{table}

\subsection{Large wave numbers}
We finally investigate the performance of the WG method for the
Helmholtz equation with large wave numbers. The homogeneous
Helmholtz problem studied in the Subsection \ref{convex} is employed
again for this purpose. Also, the $RT_0$ and $RT_1$ elements are
used to solve the homogeneous Helmholtz equation with the Robin
boundary condition. Since this problem is defined on a structured
hexagon domain, a uniform triangular mesh with a constant mesh size
$h$ throughout the domain is used. This enables us to precisely
evaluate the impact of the mesh refinements. Following the
literature works \cite{bao04,fw}, we will focus only on the relative
$H^1$ semi-norm in the present study.

\begin{figure}[!h]
\centering
\begin{tabular}{cc}
  \resizebox{2.45in}{2.15in}{\includegraphics{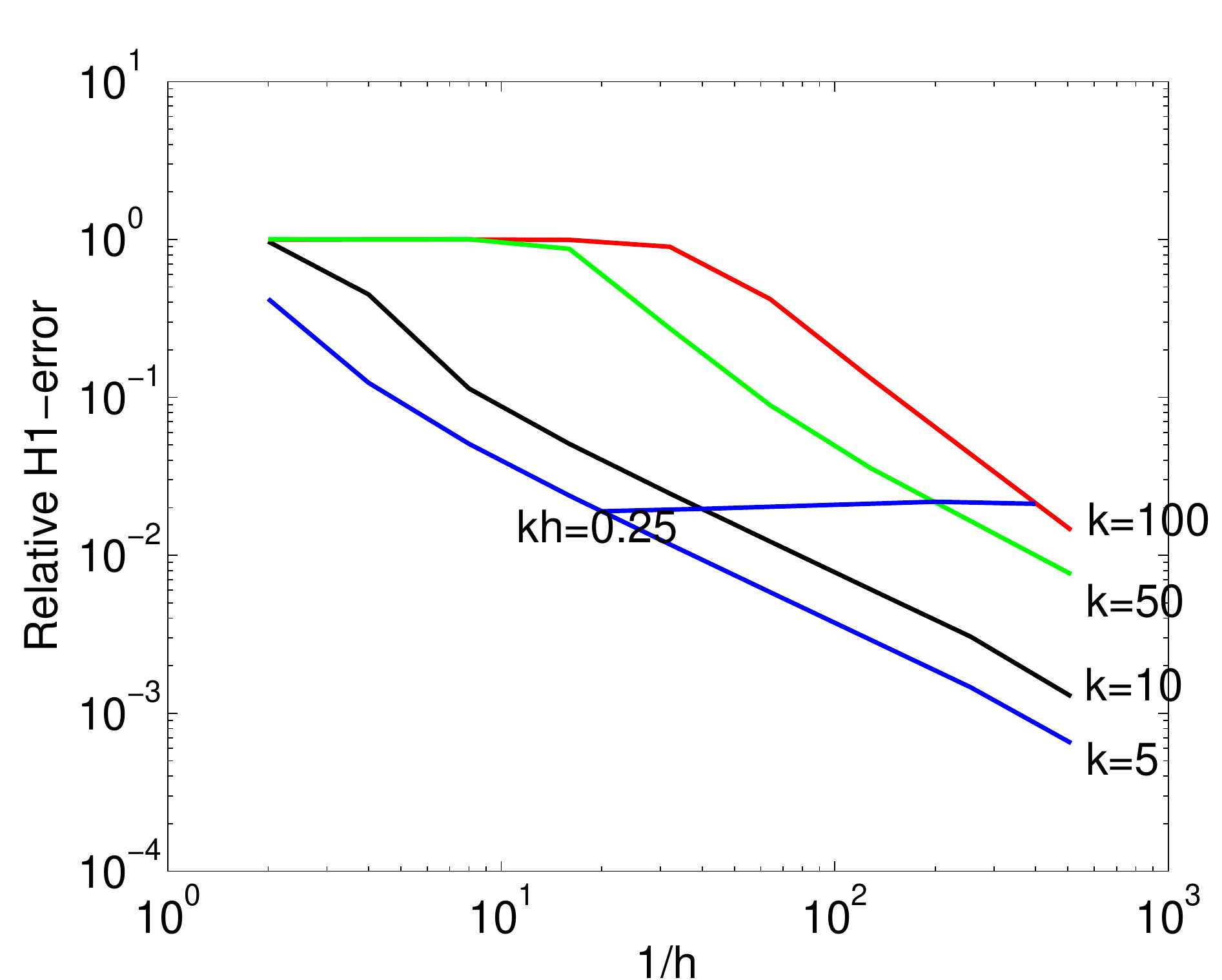}}
  \resizebox{2.45in}{2.15in}{\includegraphics{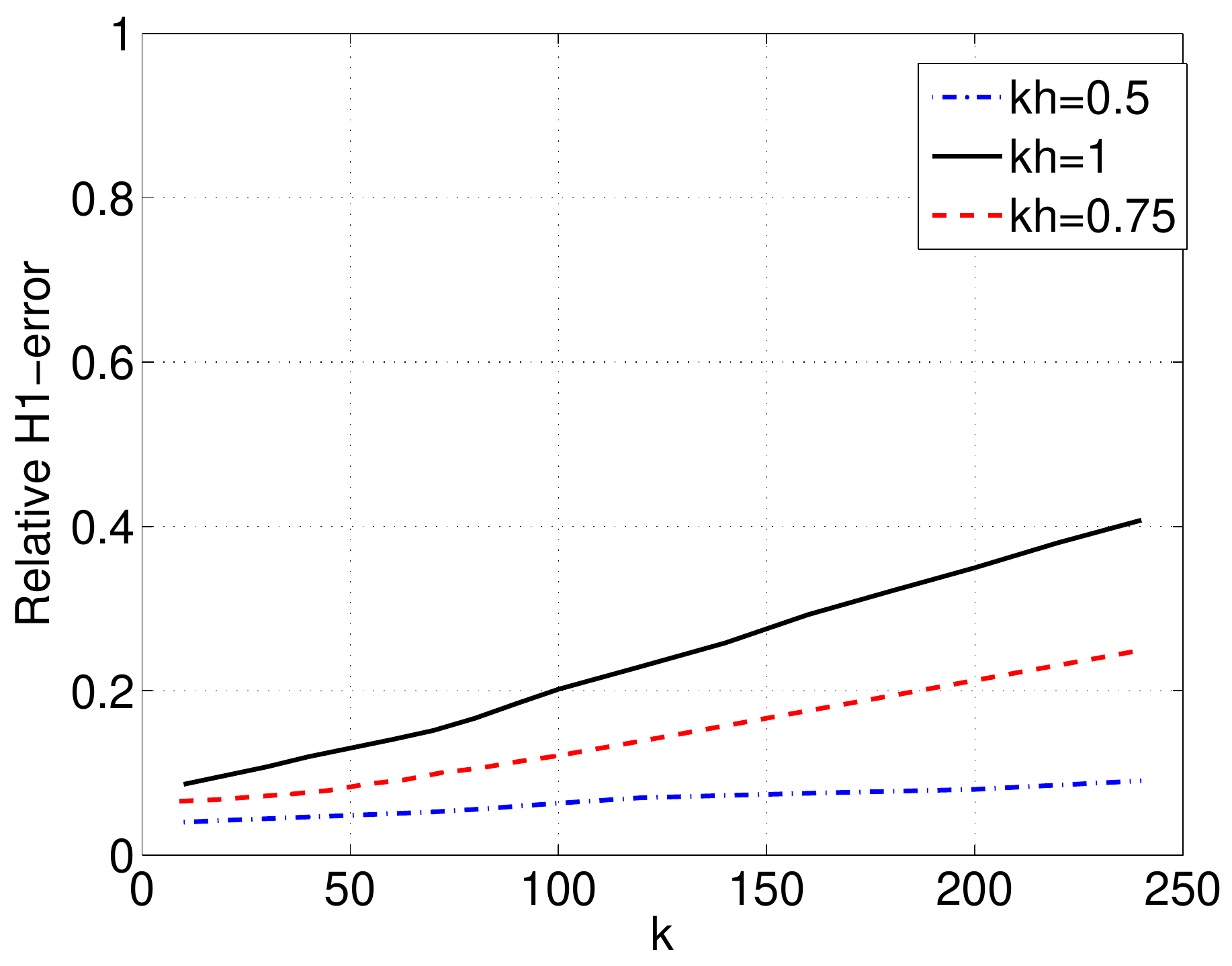}}
\end{tabular}
\caption{Relative $H^1$ error of the WG solution. Left: with respect
to $1/h$; Right: with respect to wave number $k$. } \label{fig.kh}
\end{figure}

To study the non-robustness behavior with respect to the wave number
$k$, i.e., the pollution effect, we solve the corresponding
Helmholtz equation by using piecewise constant WG method with
various mesh sizes for four wave numbers $k=5$, $k=10$, $k=50$, and
$k=100$, see Fig. \ref{fig.kh} (left) for the WG performance. From
Fig. \ref{fig.kh} (left), it can be seen that when $h$ is smaller,
the WG method immediately begins to converge for the cases $k=5$ and
$k=10$. However, for large wave numbers $k=50$ and $k=100$, the
relative error remains to be about $100$\%, until $h$ becomes to be quite
small or $1/h$ is large. This indicates the presence of the
pollution effect. In the same figure, we also show the errors of
different $k$ values by fixing $kh=0.25$. Surprisingly, we found
that the relative $H^1$ error does not evidently increase as $k$
becomes larger. The convergence line for $kh=0.25$ looks almost
flat, with a very little slope. We note that the present result of
the WG method is as good as the one reported in \cite{fw} by using a
penalized discontinuous Galerkin approach with optimized parameter
values \cite{fw}. We emphasize that WG has advantage over the
penalized DG in that no parameters are involved in the numerical
scheme.

On the other hand, the good performance of the WG method for the
case $kh=0.25$ does not mean that the WG method could be free of
pollution effect. In fact, it is known theoretically  \cite{BabSau}
that the pollution error cannot be eliminated completely in two- and
higher-dimensional spaces for Galerkin finite element methods. In
the right chart of Fig. \ref{fig.kh}, we examine the numerical
errors by increasing $k$, under the constraint that $kh$ is a
constant. Huge wave numbers, up to $k=240$, are tested. It can be
seen that when the constant changes from $0.5$ to $0.75$ and $1.0$,
the non-robustness behavior against $k$ becomes more and more
evident. However, the slopes of $kh$=constant lines remain to be
small and the increment pattern with respect to $k$ is always
monotonic. This suggests that the pollution error is well controlled
in the WG solution.

\begin{figure}[!h]
\centering
\begin{tabular}{cc}
  \resizebox{2.45in}{2.15in}{\includegraphics{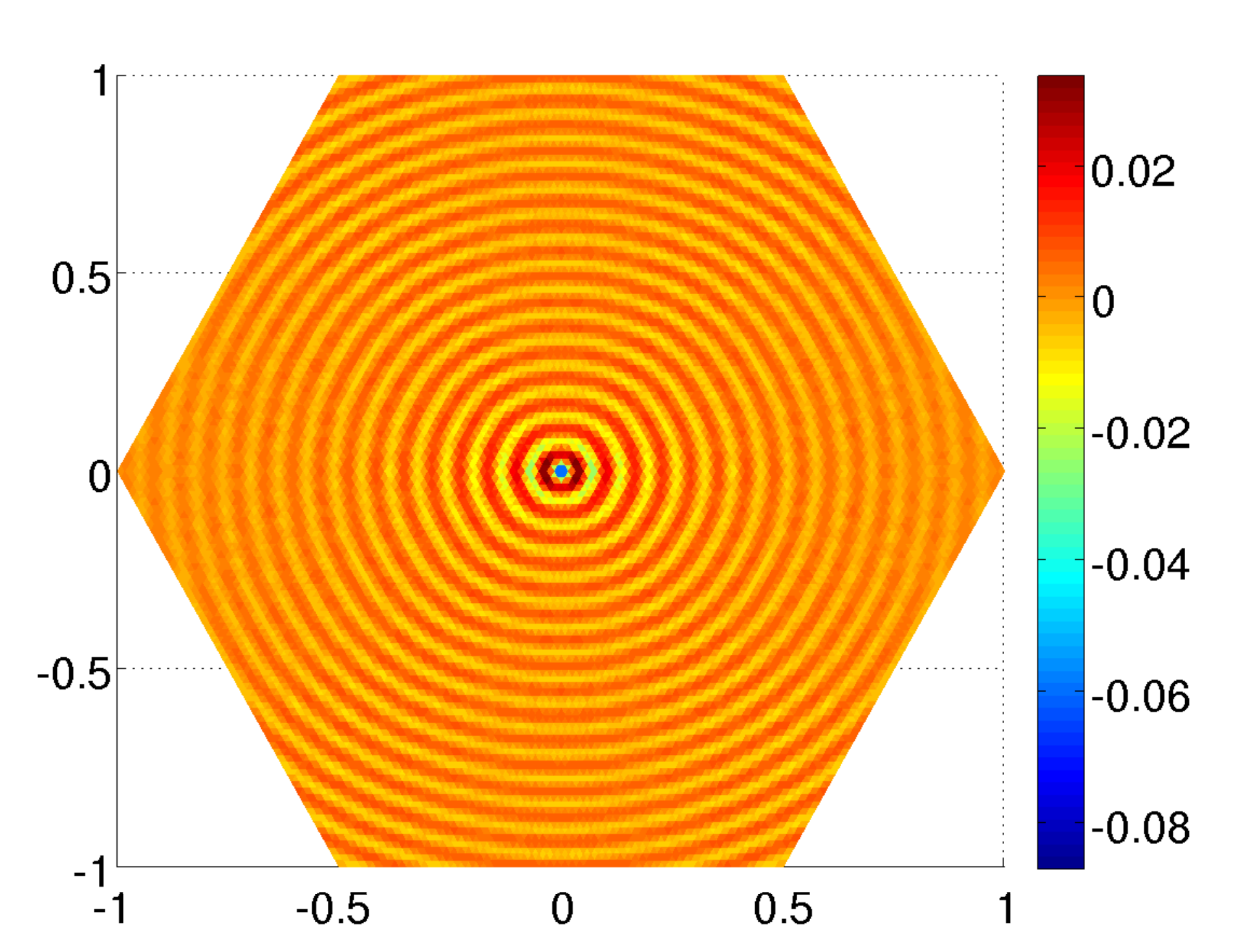}}
  \resizebox{2.45in}{2.15in}{\includegraphics{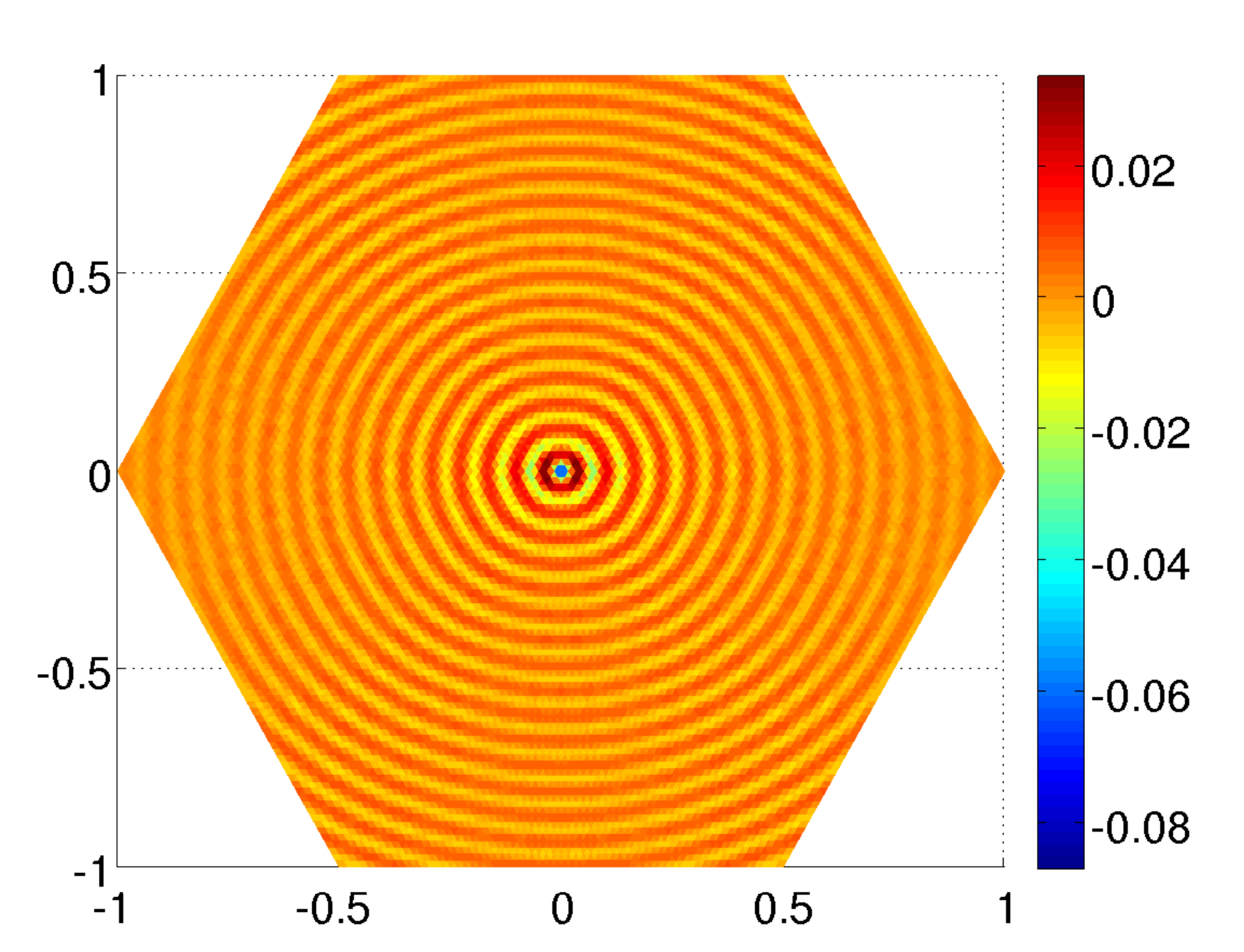}}
\end{tabular}
\caption{Exact solution (left) and piecewise constant WG
approximation (right) for $k=100,$ and $h=1/60.$} \label{fig.solu2d}
\end{figure}

\begin{figure}[!h]
\centering
\begin{tabular}{c}
  \resizebox{2.35in}{2.15in}{\includegraphics{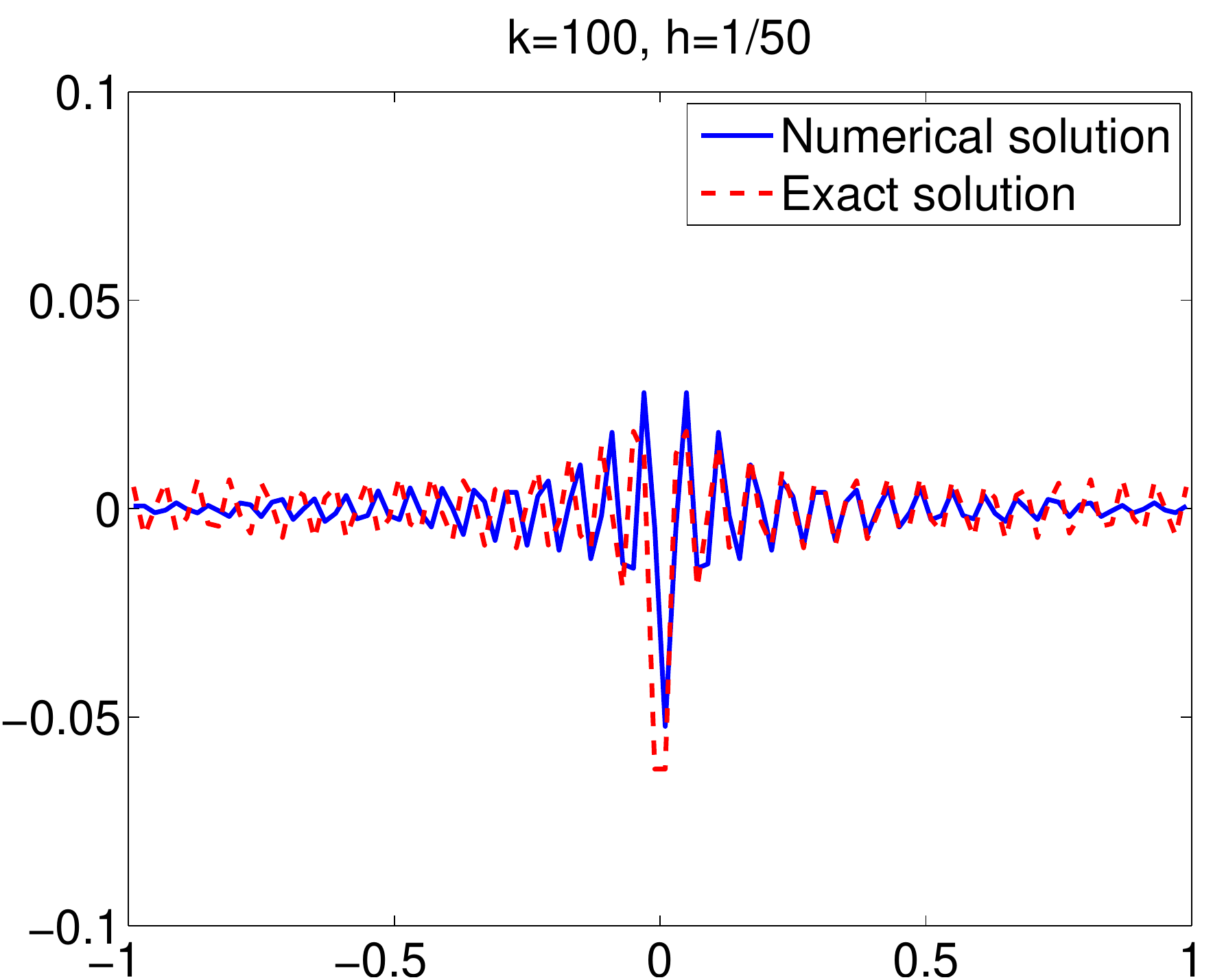}}\\
  \resizebox{2.35in}{2.15in}{\includegraphics{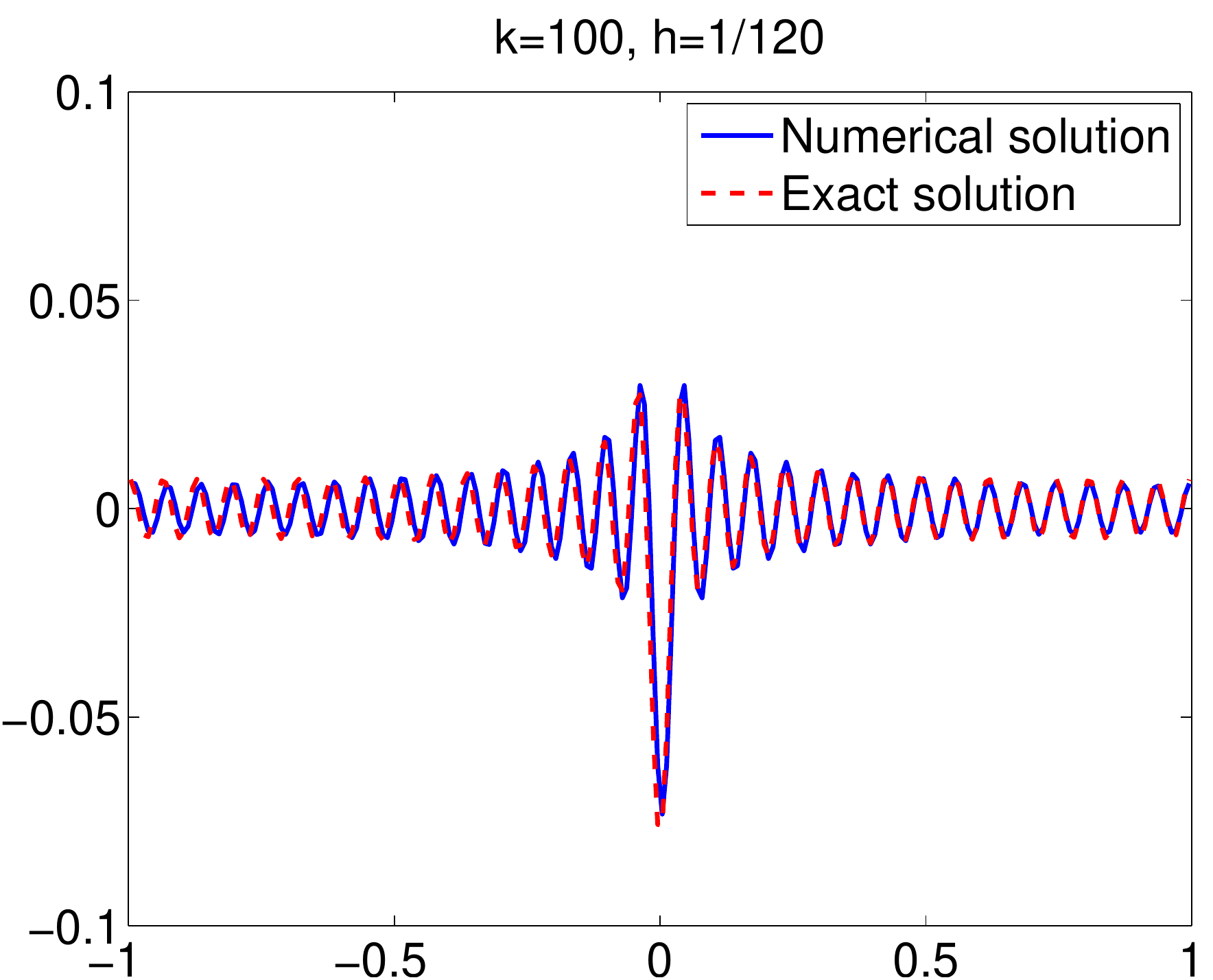}}\\
  \resizebox{2.35in}{2.15in}{\includegraphics{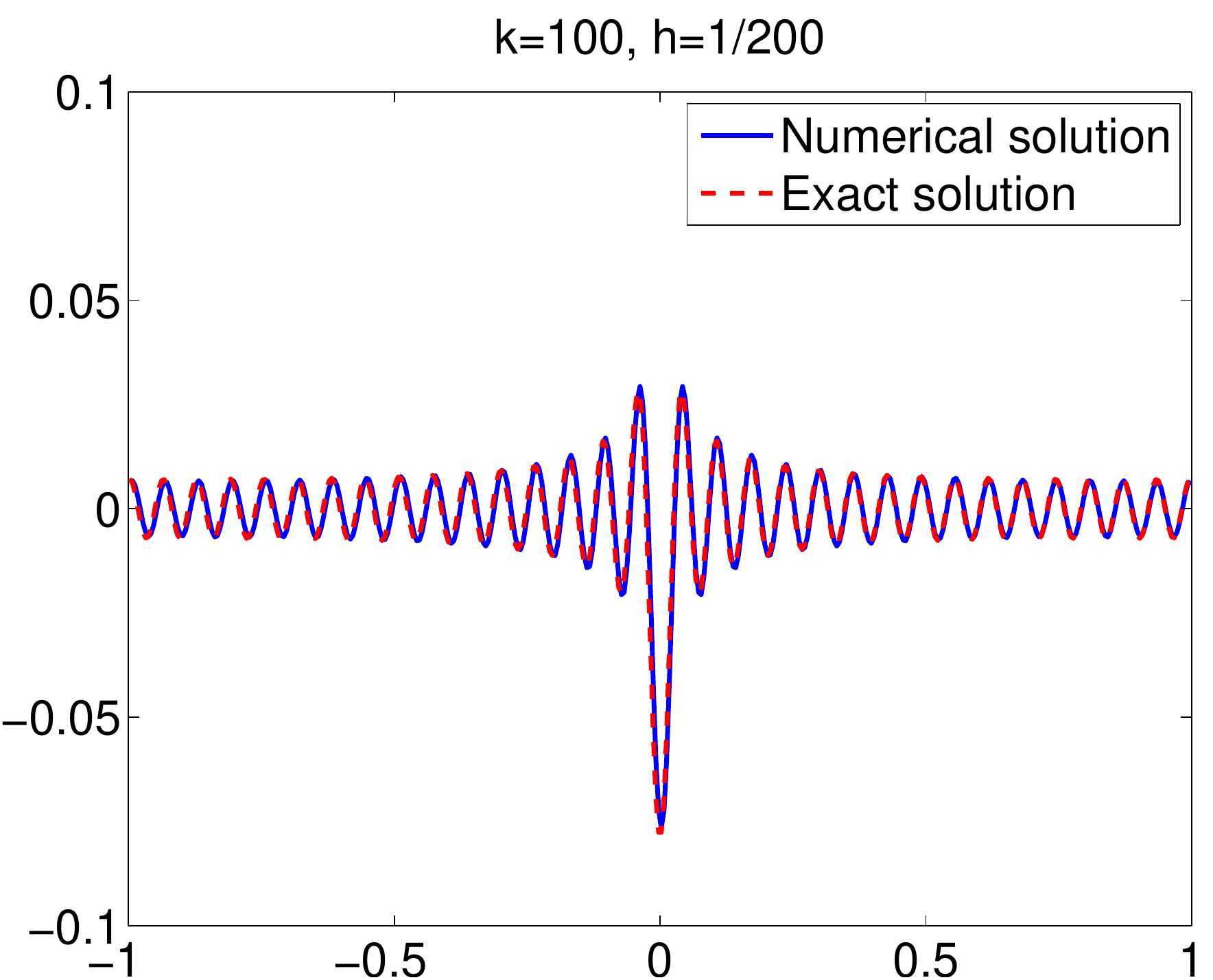}}
\end{tabular}
\caption{The trace plot along $x$-axis or $y=0$ form WG solution
using piecewise constants. } \label{fig.solu1d}
\end{figure}

\begin{figure}[!h]
\centering
\begin{tabular}{cc}
  \resizebox{2.45in}{2.15in}{\includegraphics{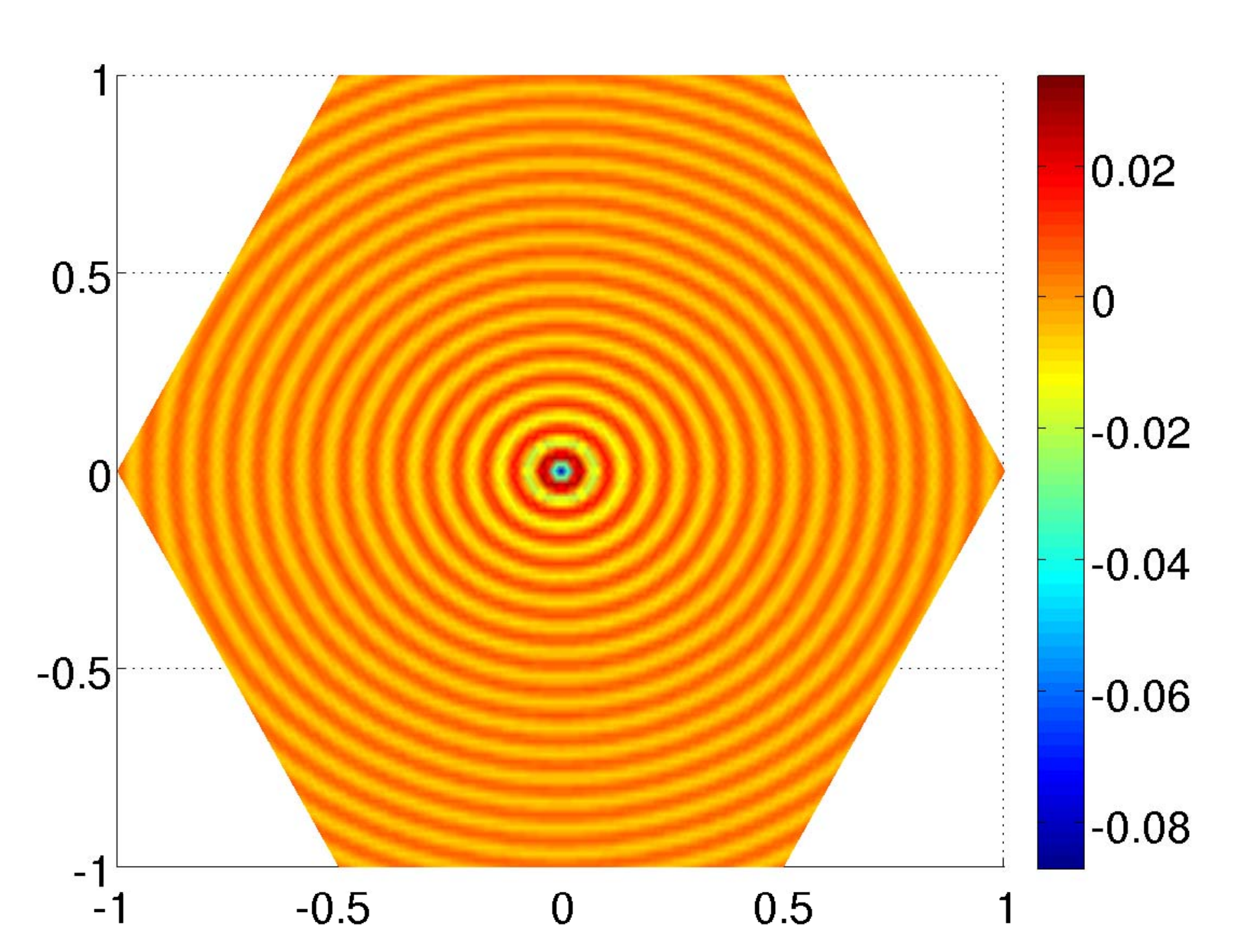}}
  \resizebox{2.45in}{2.15in}{\includegraphics{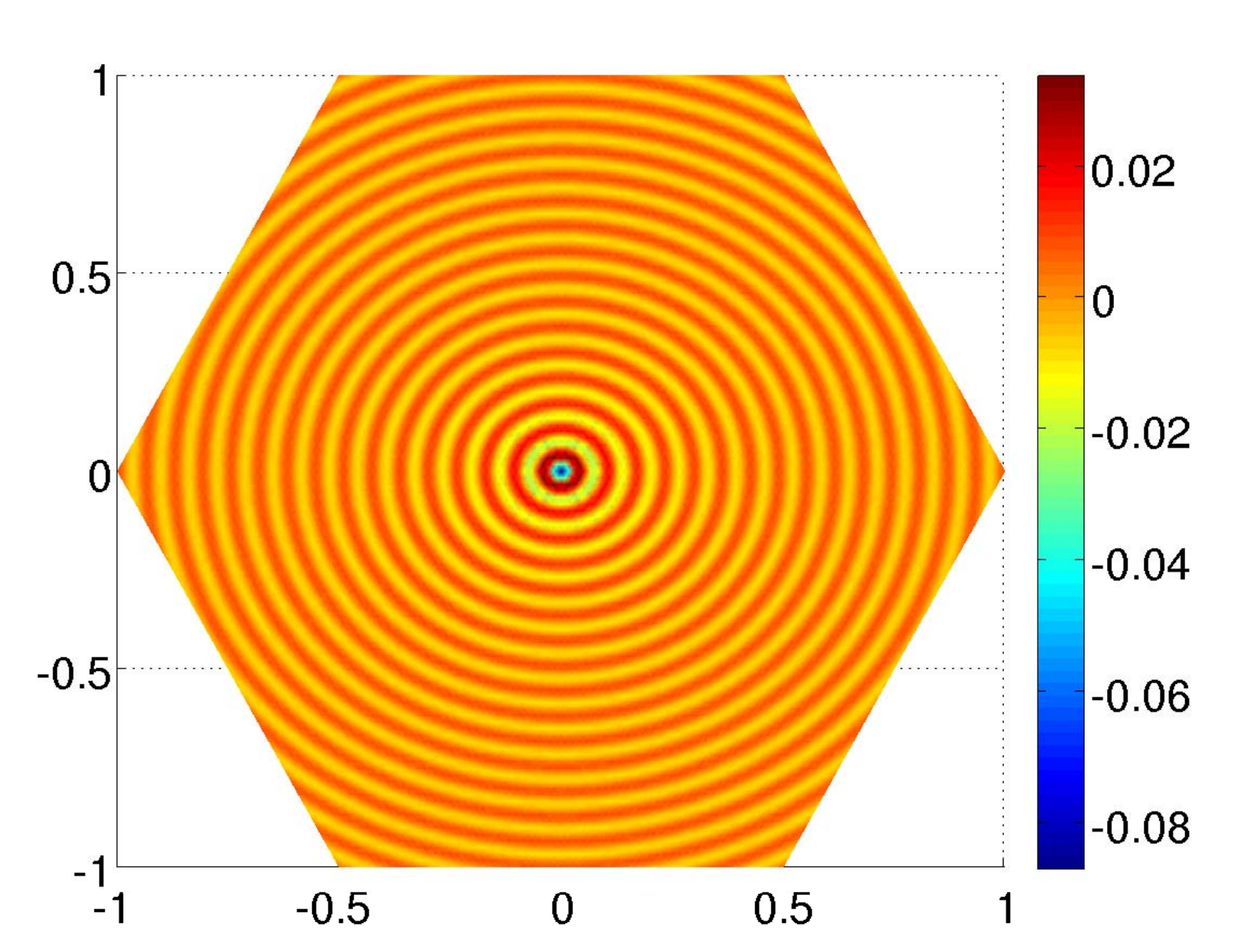}}
\end{tabular}
\caption{Exact solution (left) and piecewise linear WG approximation
(right) for $k=100,$ and $h=1/60.$} \label{fig.solu2d_linear}
\end{figure}

\begin{figure}[!h]
\centering
\begin{tabular}{c}
  \resizebox{2.35in}{2.15in}{\includegraphics{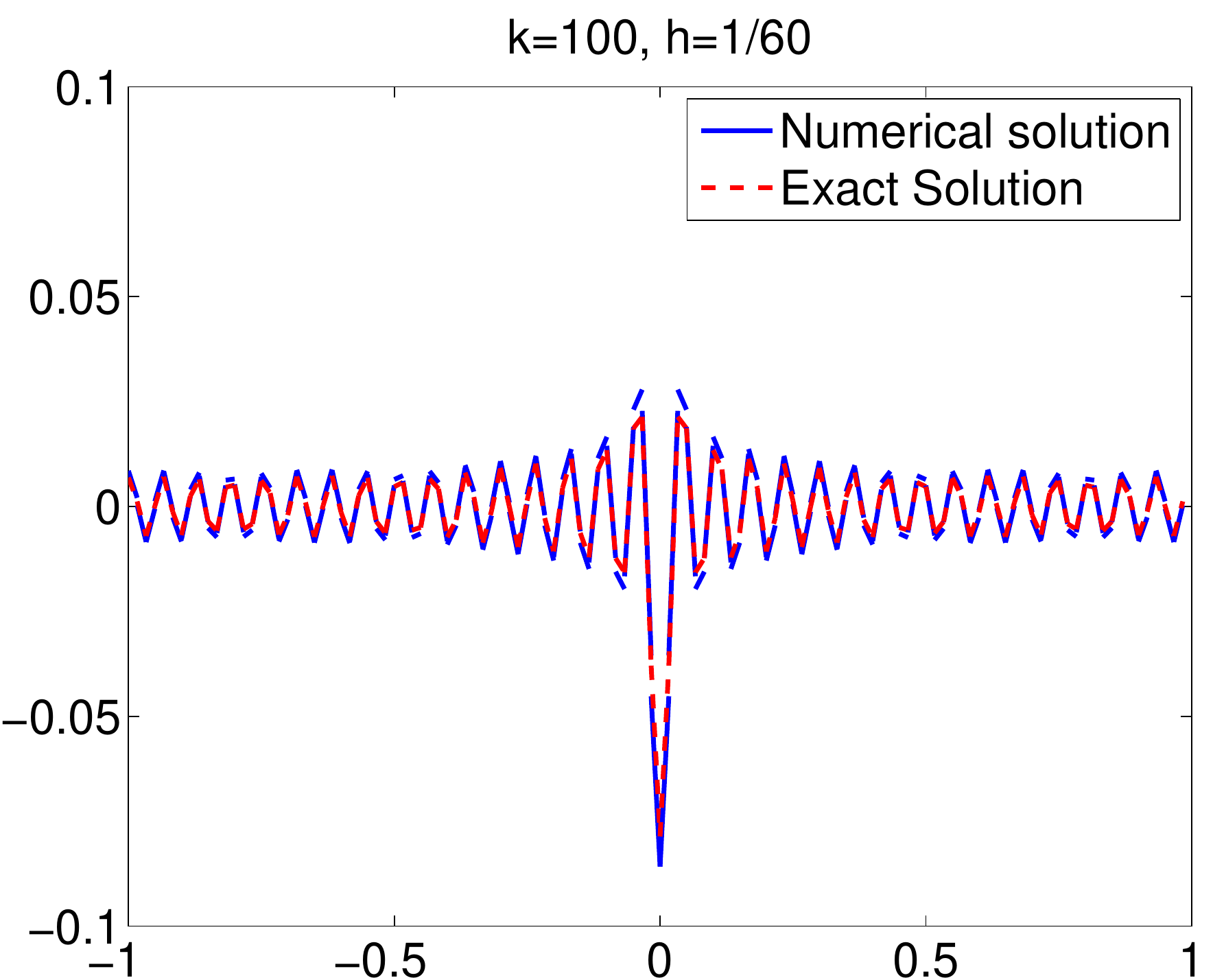}}\\
\end{tabular}
\caption{The trace plot along $x$-axis or $y=0$ form WG solution
using piecewise linear elements. } \label{fig.solu1d_linear}
\end{figure}

In the rest of the paper, we shall present some numerical results
for the WG method when applied to a challenging case of high wave
numbers. In Fig. \ref{fig.solu2d} and \ref{fig.solu2d_linear}, the
WG numerical solutions are plotted against the exact solution of the
Helmholtz problem. Here we take a wave number $k=100$ and mesh size
$h=1/60$ which is relatively a coarse mesh. With such a coarse mesh,
the WG method can still capture the fast oscillation of the
solution. However, the numerically predicted magnitude of the
oscillation is slightly damped for waves away from the center when
piecewise constant elements are employed in the WG method. Such
damping can be seen in a trace plot along $x$-axis or $y=0$. To see
this, we consider an even worse case with $k=100$ and $h=1/50$. The
result is shown in the first chart of Fig. \ref{fig.solu1d}. We note
that the numerical solution is excellent around the center of the
region, but it gets worse as one moves closer to the boundary. If we
choose a smaller mesh size $h=1/120$, the visual difference between
the exact and WG solutions becomes very small, as illustrate in Fig.
\ref{fig.solu1d}. If we further choose a mesh size $h=1/200$, the
exact solution and the WG approximation look very close to each
other. This indicates an excellent convergence of the WG method when
the mesh is refined. In addition to mesh refinement, one may also
obtain a fast convergence by using high order elements in the WG
method. Figure \ref{fig.solu1d_linear} illustrates a trace plot for
the case of $k=100$ and $h=1/60$ when piecewise linear elements are
employed in the WG method. It can be seen that the computational
result with this relatively coarse mesh captures both the fast
oscillation and the magnitude of the exact solution very well.

\section{Concluding Remarks}

The numerical experiments indicate that the WG method as introduced
in \cite{wy} is a very promising numerical technique for solving the
Helmholtz equations with large wave numbers. This finite element
method is robust, efficient, and easy to implement. On the other
hand, a theoretical investigation for the WG method should be
conducted by taking into account some useful features of the
Helmholtz equation when special test functions are used. It would
also be valuable to test the performance of the WG method when high
order finite elements are employed to the Helmholtz equations
with large wave numbers in two and three dimensional spaces.

%\begin{figure}[!h]
%\centering
%\begin{tabular}{cc}
%  \resizebox{2.45in}{2.15in}{\includegraphics{hex_errorVSk_L2error.eps}}
%  \resizebox{2.45in}{2.15in}{\includegraphics{hex_errorVSk_H1error.eps}}
%\end{tabular}
%\caption{Error profile.}
%
%\end{figure}

\newpage
\, \

\end{document}